\documentclass[12pt,oneside,reqno]{amsart}
\usepackage{mathrsfs}
\usepackage{graphics}
\usepackage{amssymb}
\usepackage{color}
\usepackage{enumerate}
\newcommand{\dif}{\mathrm{d}}

\newcommand{\be}{\begin{eqnarray}}
\newcommand{\ee}{\end{eqnarray}}
\newcommand{\ce}{\begin{eqnarray*}}
\newcommand{\de}{\end{eqnarray*}}
\newtheorem{theorem}{Theorem}[section]
\newtheorem{lemma}[theorem]{Lemma}
\newtheorem{remark}[theorem]{Remark}
\newtheorem{definition}[theorem]{Definition}
\newtheorem{proposition}[theorem]{Proposition}
\newtheorem{Example}[theorem]{Example}
\newtheorem{corollary}[theorem]{Corollary}
\newtheorem{condition}[theorem]{Condition}
\def\e{\varepsilon}
\def\s{\sigma}

\def\a{\alpha}

\def\b{\beta}
\def\d{\delta}
\def\p{\partial}
\def\g{\gamma}
\def\l{\lambda}

\def\[{{\Big[}}
\def\]{{\Big]}}
\def\<{{\langle}}
\def\>{{\rangle}}
\def\({{\Big(}}
\def\){{\Big)}}

\def\no{\nonumber}
\def\bt{\begin{theorem}}
\def\et{\end{theorem}}
\def\bl{\begin{lemma}}
\def\el{\end{lemma}}
\def\br{\begin{remark}}
\def\er{\end{remark}}
\def\bx{\begin{Example}}
\def\ex{\end{Example}}
\def\bd{\begin{definition}}
\def\ed{\end{definition}}
\def\bp{\begin{proposition}}
\def\ep{\end{proposition}}
\def\bc{\begin{corollary}}
\def\ec{\end{corollary}}
\def\bco{\begin{condition}}
\def\eco{\end{condition}}

\def\cD{{\mathcal D}}

\def\cG{{\mathcal G}}

\def\cN{{\mathcal N}}
\def\cO{{\mathcal O}}
\def\cP{{\mathcal P}}

\def\mE{{\mathbb E}}
\def\mF{{\mathbb F}}

\def\mH{{\mathbb H}}

\def\mN{{\mathbb N}}

\def\mP{{\mathbb P}}

\def\mR{{\mathbb R}}
\def\mS{{\mathbb S}}

\def\mW{{\mathbb W}}

\def\sA{{\mathscr A}}
\def\sB{{\mathscr B}}

\def\sF{{\mathscr F}}
\def\sG{{\mathscr G}}

\def\sL{{\mathscr L}}

\def\sV{{\mathscr V}}

\def\geq{\geqslant}
\def\leq{\leqslant}

\begin{document}

\allowdisplaybreaks

\title{Large deviations of multiscale multivalued McKean-Vlasov stochastic systems}

\author{Huijie Qiao}

\dedicatory{School of Mathematics,
Southeast University,\\
Nanjing, Jiangsu 211189, P.R.China\\
hjqiaogean@seu.edu.cn}



\thanks{This work was supported by NSF of China (No.12071071) and the Jiangsu Provincial Scientific Research Center of Applied Mathematics (No. BK20233002).}


\date{}

\begin{abstract}
This work concerns about multiscale multivalued McKean-Vlasov stochastic systems. First of all, we use a contractive mapping principle to establish the well-posedness for fully coupled multivalued McKean-Vlasov stochastic systems under non-Lipschitz conditions. Then for multiscale multivalued McKean-Vlasov stochastic systems with small noises, we prove a large deviation principle by a weak convergence approach. As a by-product, two averaging principles are obtained.
\end{abstract}

\maketitle \rm

\section{Introduction}

The asymptotic theory of large deviation principles (LDPs for short) quantifies the rate of convergence for the probability of rare events. And a Freidlin-Wentzell LDP provides an estimate for the probability that the sample path of an It\^o diffusion will stray far from the mean path when the size of the driving Brownian motion is small with respect to a pathspace norm. Freidlin-Wentzell LDPs have been established for many equations, such as multiscale stochastic differential equations (SDEs for short) (\cite{ds, kp, rl, aap, av2, av3}), multiscale multivalued SDEs (\cite{ku, q1}), multivalued McKean-Vlasov SDEs (\cite{arrst, flqz}) and multiscale McKean-Vlasov SDEs (\cite{bs1, bs3, ghl, hlls}). 

Now, consider the system of slow-fast multivalued McKean-Vlasov SDEs on $\mR^{n} \times \mR^{m}$: for any $T>0$
\be\left\{\begin{array}{l}
\dif X_{t}^{\e,\d}\in -A_1(X_{t}^{\e,\d})\dif t+b_{1}(X_{t}^{\e,\d},\sL_{X_{t}^{\e,\d}},Y_{t}^{\e,\d})\dif t+\sqrt\e\s_{1}(X_{t}^{\e,\d},\sL_{X_{t}^{\e,\d}},Y_{t}^{\e,\d})\dif W^1_{t},\\
X_{0}^{\e,\d}=x_0\in\overline{\cD(A_1)},\quad  0\leq t\leq T,\\
\dif Y_{t}^{\e,\d}\in -A_2(Y_{t}^{\e,\d})\dif t+\frac{1}{\d}b_{2}(X_{t}^{\e,\d},\sL_{X_{t}^{\e,\d}},Y_{t}^{\e,\d})\dif t+\frac{1}{\sqrt{\d}}\s_{2}(X_{t}^{\e,\d},\sL_{X_{t}^{\e,\d}},Y_{t}^{\e,\d})\dif W^2_{t},\\
Y_{0}^{\e,\d}=y_0\in\overline{\cD(A_2)},\quad  0\leq t\leq T.
\end{array}
\right.
\label{Eqin1}
\ee
The system (\ref{Eqin1}) is defined on a filtered probability space $(\Omega,\sF,\{\sF_{t}\}_{t\in[0,T]},\mP)$ and $(W^1_{t}), (W^2_{t})$ are $d_1$- and $d_2$-dimensional standard Brownian motions defined on it, respectively. $A_1, A_2$ are two maximal monotone operators (cf. Subsection \ref{mmo}), these mappings $b_{1}:\mR^{n}\times\cP_2(\mR^n)\times\mR^{m}\rightarrow\mR^{n}$, $\s_{1} :\mR^{n}\times\cP_2(\mR^n)\times\mR^{m}\rightarrow\mR^{n\times d_1}$, $b_{2} :\mR^{n}\times\cP_2(\mR^n)\times\mR^{m}\rightarrow\mR^{m}$, $\s_{2} :\mR^{n}\times\cP_2(\mR^n)\times\mR^{m}\rightarrow\mR^{m\times d_2}$ are all Borel measurable, and $\cP_{2}(\mR^n)$ is the set of probability measures on $\sB(\mR^n)$ with finite second moments (cf. Subsection \ref{nn}). $\sL_{X^{\e,\d}_{t}}$ denotes the distribution of $X^{\e,\d}_{t}$ under the probability measure $\mP$. Here $0<\e<1$ is a small parameter, and $\d$ represents the other small parameter which characterizes the ratio of timescales between processes $X_{\cdot}^{\e,\d}$ and $Y_{\cdot}^{\e,\d}$. 

These systems like (\ref{Eqin1}) have appeared as models for many fields, such as the control theory, biology and chemistry (cf. \cite{ahlw, lw}). Moreover, for the system (\ref{Eqin1}), to the best knowledge of the author, there are not related LDP results. Therefore, we study the LDP for the system (\ref{Eqin1}). However, in the process of achieving the goal there are two difficulties. The first difficulty lies in that those direct methods used in \cite{ghl, hlls} don't work, since the system (\ref{Eqin1}) contains two maximal monotone operators $A_1, A_2$ and they are multivalued, nonlinear and not smooth. With the help of the It\^o formula and some properties of maximal monotone operators, we overcome this difficulty. The second difficulty is that we can not obtain some estimates, such as $\mE|X_{t}^{\e,\d}-X_{s}^{\e,\d}|^{2}$, which is easy for one-valued SDEs. Hence we apply some properties of maximal monotone operators to give some limits which are sufficient for us. Besides, it is worthwhile mentioning that the distribution $\sL_{X^{\e,\d}_{t}}$ in the coefficients of the system (\ref{Eqin1}) results in that the controlled equation depends on the distribution $\sL_{X^{\e,\d}_{t}}$ but not the distribution of the controlled process $X_{t}^{\e,\d,u}$ itself (See Remark \ref{chmv}). Thus, we must make a lot of estimates about both $X^{\e,\d}$ and $X^{\e,\d,u}$ to justify Condition \ref{cond} $(ii)$. Therefore, the proof of our result is more complicated and skilled than that of some known results (cf. \cite[Theorem 3.3]{q1}).

The novelty of this paper lies in three folds. The first fold is that we use a contractive mapping principle to prove the well-posedness for fully coupled multivalued McKean-Vlasov stochastic systems under non-Lipschitz conditions. And our proof is more concise than that in some known results (cf. \cite{cw, rsx}). The second fold is that our conditions are more general than that in some known results (cf. \cite{ghl, hlls, q1}). These conditions together with the distributions in drift coefficients and diffusion coefficients make some estimates complicated and difficult. The third fold is that as a by-product, we establish two averaging principles and present the order of convergence in a special case. Comparing our result with some known results, we find that maximal monotone operators $A_1, A_2$ reduce the order of convergence (cf. Remark \ref{conrat1}). 

The rest of this paper is organized as follows. In Section \ref{pre} we introduce some notations and concepts. And the formulation of main results is placed in Section \ref{main}. In the next three sections, we prove them. In Section \ref{exam}, we give an example to explain our results. Finally, the proofs of two inequalities are attached in Section \ref{app}.

The following convention will be used throughout the paper: $C$ with or without indices will denote different positive constants whose values may change from one place to another.

\section{Preliminary}\label{pre}

In this section, we will recall some notations and concepts.

\subsection{Notations}\label{nn}
In this subsection, we introduce some notations used in the sequel.

Let $|\cdot|, \|\cdot\|$ be the norms of a vector and a matrix, respectively. Let $\langle\cdot,\cdot\rangle$ be the inner product of vectors on $\mR^n$. $U^{*}$ denotes the transpose of the matrix $U$.

Let $C(\mR^n)$ be the set of all  functions which are continuous on $\mR^n$. $C^{1}(\mR^n)$ represents the collection of all functions in $C(\mR^n)$ with continuous derivatives of order $1$. 

Let $\sB(\mR^n)$ be the Borel $\sigma$-algebra on $\mR^n$ and $\cP({\mR^n})$ be the space of all probability measures defined on $\sB(\mR^n)$ carrying the usual topology of the weak convergence. Let $\cP_{2}(\mR^n)$ be the set of probability measures on $\sB(\mR^n)$ with finite second order moments, i.e.
$$
\cP_2\left( \mathbb{R}^n \right) :=\left\{ \mu \in \cP\left( \mathbb{R}^n \right): \mu(|\cdot|^{2}):=\int_{\mathbb{R}^n}{\left| x \right|^2\mu \left( \dif x \right) <\infty} \right\}.
$$
It is known that $\cP_2(\mR^n)$ is a Polish space endowed with the $L^2$-Wasserstein distance defined by
$$
\mathbb{W}_2(\mu,\nu):= \inf\limits_{\pi\in\Phi(\mu,\nu)}\left(\int_{\mathbb{R}^n\times\mathbb{R}^n}|x-y|^{2}\pi(\dif x,\dif y)\right)^{\frac{1}{2}}, \quad \mu , \nu\in \cP_2(\mR^n),
$$
where $\Phi(\mu,\nu)$ is the set of all couplings $\pi$ with marginal distributions $\mu$ and $\nu$. Moreover, if $\xi,\zeta$ are two random variables with distributions $\sL_\xi, \sL_\zeta$ under $\mP$, respectively,
$$
\mathbb{W}^2_2(\sL_\xi, \sL_\zeta)\leq\mE|\xi-\zeta|^2,
$$
where $\mE$ stands for the expectation with respect to $\mP$.

\subsection{Maximal monotone operators}\label{mmo}

In this subsection, we introduce maximal monotone operators.

For a multivalued operator $A: \mR^n\mapsto 2^{\mR^n}$, where $2^{\mR^n}$ stands for all the subsets of $\mR^n$, set
\ce
&&\cD(A):= \left\{x\in \mR^n: A(x) \ne \emptyset\right\},\\
&&Gr(A):= \left\{(x,y)\in \mR^{2n}:x \in \cD(A), ~ y\in A(x)\right\}.
\de
We say that $A$ is monotone if $\langle x_1 - x_2, y_1 - y_2 \rangle \geq 0$ for any $(x_1,y_1), (x_2,y_2) \in Gr(A)$, and $A$ is maximal monotone if 
$$
(x_1,y_1) \in Gr(A) \iff \langle x_1-x_2, y_1 -y_2 \rangle \geq 0, \quad \forall (x_2,y_2) \in Gr(A).
$$

We give an example to explain maximal monotone operators.

\bx\label{exmmo2}
For a closed convex subset $\mathcal{O}$ of $\mathbb{R}^n$, we suppose $\operatorname{Int}(\mathcal{O})\neq\emptyset$, where $\operatorname{Int}(\cO)$ is the interior of $\cO$. Define the indicator function of $\mathcal{O}$ as follows:
$$
I_{\mathcal{O}}(x):= \begin{cases}0, & \text { if } x \in \mathcal{O}, \\ 
+\infty, & \text { if } x \notin \mathcal{O}.\end{cases}
$$
The subdifferential operator of $I_{\mathcal{O}}$ is given by
$$
\begin{aligned}
\partial I_{\mathcal{O}}(x) & :=\left\{y \in \mathbb{R}^n:\langle y, x-z\rangle \geq 0, \forall z \in \mathcal{O}\right\} \\
& = \begin{cases}\emptyset, & \text { if } x \notin \mathcal{O}, \\
\{0\}, & \text { if } x \in \operatorname{Int}(\mathcal{O}), \\
\Lambda_x, & \text { if } x \in \partial \mathcal{O},\end{cases}
\end{aligned}
$$
where $\Lambda_x$ is the exterior normal cone at $x$. By simple deduction, we know that $\partial I_{\mathcal{O}}$ is a maximal monotone operator.
\ex
In the following, we recall some properties of a maximal monotone operator $A$ (cf.\cite{cepa1}):
\begin{enumerate}[(i)]
\item
${\rm Int}(\cD(A))$ and $\overline{\mathrm{\cD}(A)}$ are convex subsets of $\mR^n$ with ${\rm Int}\left( \overline{\mathrm{\cD}(A)} \right) = {\rm Int}\( \mathrm{\cD}(A) \) 
$, where ${\rm Int}(\cD(A))$ denotes the interior of the set $\cD(A)$. 
\item For every $x\in\mR^n$, $A(x)$ is a closed and convex subset of $\mR^n$.
\end{enumerate}

Take any $T>0$ and fix it. Let $\sV_{0}$ be the set of all continuous functions $K: [0,T]\mapsto\mR^n$ with finite variations and $K_{0} = 0$. For $K\in\sV_0$ and $s\in [0,T]$, we shall use $|K|_{0}^{s}$ to denote the variation of $K$ on $[0,s]$
and write $|K|_{TV}:=|K|_{0}^{T}$. Set
\ce
&&\sA:=\Big\{(X,K): X\in C([0,T],\overline{\cD(A)}), K \in \sV_0, \\
&&\qquad\qquad\quad~\mbox{and}~\langle X_{t}-x, \dif K_{t}-y\dif t\rangle \geq 0 ~\mbox{for any}~ (x,y)\in Gr(A)\Big\}.
\de
And about $\sA$ we have two following results (cf.\cite{cepa2, ZXCH}).

\bl\label{equi}
For $X\in C([0,T],\overline{\cD(A)})$ and $K\in \sV_{0}$, the following statements are equivalent:
\begin{enumerate}[(i)]
	\item $(X,K)\in \sA$.
	\item For any $(x,y)\in C([0,T],\mR^n)$ with $(x_t, y_t)\in Gr(A)$, it holds that 
	$$
	\left\langle X_t-x_t, \dif K_t-y_t\dif t\right\rangle \geq0.
	$$
	\item For any $(X^{'},K^{'})\in \sA$, it holds that 
	$$
	\left\langle X_t-X_t^{'},\dif K_t-\dif K_t^{'}\right\rangle \geq0.
	$$
\end{enumerate}
\el

\bl\label{inteineq}
Assume that $\text{Int}(\cD(A))\ne\emptyset$. For any $a\in \text{Int}(\cD(A))$, there exist $M_1 >0$, and $M_{2},M_{3}\geq0$ such that  for any $(X,K)\in \sA$ and $0\leq s<t\leq T$,
$$
\int_s^t{\left< X_r-a, \dif K_r \right>}\geq M_1\left| K \right|_{s}^{t}-M_2\int_s^t{\left| X_r-a\right|}\dif r-M_3\left( t-s \right) .
$$
\el

\subsection{Multivalued McKean-Vlasov SDEs}

In this subsection, we introduce multivalued McKean-Vlasov SDEs. 

Fix $T>0$ and consider the following multivalued McKean-Vlasov SDE on $\mR^n$:
\be
\dif X_t\in -A(X_t)\dif t+b(X_t,\sL_{X_t})\dif t+\s(X_t,\sL_{X_t})\dif W^1_t, \quad 0\leq t\leq T, \label{eq1}
\ee
where $A$ is a maximal monotone operator with $\text{Int}(\cD(A))\ne\emptyset$, the coefficients $b: \mR^n\times\cP_2(\mR^n)\mapsto{\mR^n}, \,\,\sigma:\mR^n\times\cP_2(\mR^n)\mapsto{\mR^n}\times{\mR^{d_1}}$ are Borel measurable and $W^1_{\cdot}$ is a $d_1$-dimensional Brownian motion on a filtered probability space $(\Omega, \mathscr{F}, \{\mathscr{F}_t\}_{t\in[0,T]}, \mP)$.

\bd\label{strosolu}
We say that Eq.$(\ref{eq1})$ admits a strong solution with the initial value $X_0\in\overline{\cD(A)}$ if there exists a pair of adapted processes $(X,K)$ on $(\Omega, \mathscr{F}, \{\mathscr{F}_t\}_{t\in[0,T]}, \mP)$ such that

(i) $X_t\in{\mathscr{F}_t^{W^1}}$, where $\{\mathscr{F}_t^{W^1}\}_{t\in[0,T]}$ stands for the $\sigma$-field filtration generated by $W^1$,

(ii) $(X_{\cdot}(\omega),K_{\cdot}(\omega))\in \sA$ a.s. $\mP$,

(iii) it holds that
\ce
\mP\left\{\int_0^T(\mid{b(X_s,\sL_{X_s})}\mid+\parallel{\sigma(X_s,\sL_{X_s})}\parallel^2)\dif s<+\infty\right\}=1,
\de
and
\ce
X_t=X_0-K_{t}+\int_0^tb(X_s,\sL_{X_s})\dif s+\int_0^t\sigma(X_s,\sL_{X_s})\dif W^1_s, \quad 0\leq{t}\leq{T}, \quad a.s.~\mP.
\de
\ed

\subsection{A general criterion of large deviation principles}

In this subsection, we present a general criterion to establish the large deviation principle. 

Let $(\mS,\rho)$ be a Polish space. For each $\e>0$, let $X^{\e}$ be a $\mS$-valued random variable given on $(\Omega, \mathscr{F}, \{\mathscr{F}_t\}_{t\in[0,T]}, \mP)$.

\bd\label{compleve}
The function $I: \mathbb{S}\mapsto [0,\infty]$ is called a rate function if $I$ is lower semicontinuous. Moreover, a rate function $I$ is called a good rate function if for each $M<\infty$, $\{\varsigma\in \mathbb{S}:I(\varsigma)\leq M\}$ is a compact subset of $\mathbb{S}$.
\ed

\bd
We say that $\{X^{\e}\}$ satisfies the large deviation principle with the speed $\e^{-1}$ and the good rate function $I$, if for any subset $B\in \sB(\mS)$,
$$
-\inf\limits_{\varsigma\in {\rm Int}(B)}I(\varsigma)\leq\liminf_{\e\rightarrow 0}\e\log\mP(X^{\e}\in{\rm Int}(B))\leq \limsup\limits_{\e\rightarrow 0}\e\log\mP(X^{\e}\in \bar{B})\leq -\inf\limits_{\varsigma\in \bar{B}}I(\varsigma),
$$
where ${\rm Int}(B)$ and $\bar{B}$ denote the interior and the closure of $B$, respectively and they are taken in $\mS$.
\ed

\bd
We say that $\{X^{\e}\}$ satisfies the Laplace principle with the speed $\e^{-1}$ and the good rate function $I$,  if for any real bounded continuous function $G$ on $\mathbb{S}$,
\ce
\lim\limits_{\varepsilon\rightarrow 0}\varepsilon \log \mE\left\{\exp\left[-\frac{G(X^{\e})}{\e}\right]\right\}=-\inf\limits_{\varsigma\in \mathbb{S}}\(G(\varsigma)+I(\varsigma)\).
\de
\ed

In the sequel, we take $\mS:=C([0,T],\overline{\cD(A_1)})$ and $\rho(X^1, X^2):=\sup\limits_{t\in[0,T]}|X^1_t-X^2_t|$ for $X^1, X^2\in C([0,T],\overline{\cD(A_1)})$. Note that the large deviation principle is equivalent to the Laplace principle (cf. \cite{de}). Therefore, in order to obtain the large deviation principle for $\{X^{\e}\}$, we prove the Laplace principle for $\{X^{\e}\}$. Then we state the conditions under which the Laplace principle holds. Let $\mH$ be the Hilbert space of absolutely continuous functions from $[0,T]$ to $\mR^{d_1+d_2}$ with square integrable derivatives, i.e.
$$
\mathbb{H}:=\left\{h: [0,T]\rightarrow\mR^{d_1+d_2}; \|h\|^2_{\mathbb{H}}:=\int_{0}^{T}|\dot{h}(t)|^{2}\dif t\right\},
$$ 
where $\dot{h}$ stands for the generalized derivative of $h$. Let $\mathcal{A}$ be the collection of predictable processes $u(\omega, \cdot)$ belonging to $\mathbb{H}$ a.s. $\omega$. For each $N\in\mN$ we define two following spaces:
\ce
\mathbf{D}_{2}^{N}:=\left\{h\in\mathbb{H}: \|h\|_{\mathbb{H}}^{2}\leq N \right\}, \quad \mathbf{A}_{2}^{N}:=\left\{u\in\mathcal{A}: u(\omega, \cdot)\in\mathbf{D}_{2}^{N}, a.s.~\omega \right\}.
\de
We equip $\mathbf{D}_{2}^{N}$ with the weak convergence topology in $\mH$. So, $\mathbf{D}_{2}^{N}$ is metrizable as a compact Polish space. In the sequel, $\mathbf{D}_{2}^{N}$ will be always endowed with this topology.

\bco\label{cond}
Let $\cG^{\e} : C([0,T];\mathbb{R}^{d_1+d_2})\mapsto\mS$ be  a family of measurable mappings. There exists a
measurable mapping $\cG^{0} : \mH\mapsto\mS$ such that

$(i)$ for $\{h_{\e}, \e>0\}\subset\mathbf{D}_2^{N}$, $h\in \mathbf{D}_2^{N}$, if $h_{\e}\rightarrow h$ as $\e\rightarrow 0$, then
\ce
\cG^{0}\left(h_{\e}\right)\longrightarrow \cG^{0}\left(h\right).
\de

$(ii)$ for $\{u_{\e},\e>0\}\subset \mathbf{A}_{2}^{N}$, and any $\eta>0$,
$$
\lim\limits_{\e\rightarrow0}\mP\left(\sup\limits_{t\in[0,T]}\left|\cG^{\e}\left(\sqrt{\e}W_\cdot+u_{\e}\right)(t)-\cG^{0}\left(u_\e\right)(t)\right|>\eta\right)=0,
$$
where $W_\cdot$ is a $d_1+d_2$-dimensional Brownian motion.

\eco

 Given $\varsigma\in\mS$, let ${\bf D}_{\varsigma}=\{h\in\mH: \varsigma=\cG^{0}(h)\}$. Let $I:\mathbb{S}\mapsto [0,\infty]$ be defined by
$$
I(\varsigma)=\frac{1}{2}\inf\limits_{h\in{\bf D}_{\varsigma}}\|h\|_{\mathbb{H}}^{2}.
$$
The following result is due to \cite[Theorem 3.2]{msz}.

\bt\label{ldpbase}
Set $X^{\e}:=\cG^{\e}(\sqrt{\e}W)$. Assume that Condition \ref{cond} holds. Then $\{X^{\e}\}$ satisfies the Laplace principle with the rate function $I$ given above. In particular, $\{X^{\e}\}$ satisfies the large deviation principle with the same rate function $I$.
\et

\section{Main results}\label{main}

In this section, we formulate the main results in this paper.

\subsection{Well-posedness of multivalued McKean-Vlasov stochastic systems}

In this subsection, we present the well-posedness result for multivalued McKean-Vlasov stochastic systems.

Consider the following system on $\mR^{n} \times \mR^{m}$:
\be\left\{\begin{array}{l}
\dif X_{t}\in -A_1(X_{t})\dif t+b_{1}(X_{t},\sL_{X_{t}},Y_{t})\dif t+\s_{1}(X_{t},\sL_{X_{t}},Y_{t})\dif W^1_{t},\\
X_{0}=\xi\in\overline{\cD(A_1)},\quad  0\leq t\leq T,\\
\dif Y_{t}\in -A_2(Y_{t})\dif t+b_{2}(X_{t},\sL_{X_{t}},Y_{t})\dif t+\s_{2}(X_{t},\sL_{X_{t}},Y_{t})\dif W^2_{t},\\
Y_{0}=\zeta\in\overline{\cD(A_2)},\quad  0\leq t\leq T,
\end{array}
\right.
\label{Eqeu}
\ee
where $\xi, \zeta$ are $\sF_0$-measurable $\overline{\cD(A_1)}$- and $\overline{\cD(A_2)}$-valued random variables with $\mE|\xi|^2<\infty, \mE|\zeta|^2<\infty$, respectively, and independent of $W^1,W^2$.

Assume:
\begin{enumerate}[$(\mathbf{H}_{A_1})$]
\item $\text{Int}(\cD(A_1))\ne\emptyset$.
\end{enumerate}
\begin{enumerate}[$(\mathbf{H}^1_{b_{1}, \s_{1}})$]
\item $(i)$ For $y\in\mR^m$, $b_1(x,\mu,y)$ is continuous in $(x,\mu)$, and there exists a constant $L_{b_{1}, \s_{1}}>0$ such that for $x_1,x_2\in\mR^n$, $\mu_1, \mu_2\in\cP_{2}(\mR^n)$, $y\in\mR^m$, 
\ce
&&2\<x_1-x_2, b_{1}(x_{1},\mu_1,y)-b_{1}(x_{2},\mu_2,y)\>\leq L_{b_{1},\s_{1}}\(|x_{1}-x_{2}|^2+\mW^2_2(\mu_1,\mu_2)\),\\
&&\|\s_{1}(x_{1},\mu_1,y)-\s_{1}(x_{2},\mu_2,y)\|^2\leq L_{b_{1},\s_{1}}\(|x_{1}-x_{2}|^2+\mW^2_2(\mu_1,\mu_2)\),
\de
and for $x\in\mR^n$, $\mu\in\cP_{2}(\mR^n)$, $y_i\in\mR^m, i=1,2$
\ce
|b_{1}(x,\mu,y_1)-b_{1}(x,\mu,y_2)|^2+\|\s_{1}(x,\mu,y_1)-\s_{1}(x,\mu,y_2)\|^2\leq L_{b_{1},\s_{1}}|y_1-y_2|^2.
\de
$(ii)$ There exists a constant $\bar{L}_{b_{1}, \s_{1}}>0$ such that for $x\in\mR^n$, $\mu\in\cP_{2}(\mR^n)$, $y\in\mR^m$,
\be
|b_{1}(x,\mu,y)|^{2}+\|\s_{1}(x,\mu,y)\|^{2}\leq \bar{L}_{b_{1}, \s_{1}}(1+|x|^{2}+\mu(|\cdot|^2)+|y|^{2}).
\label{b1line}
\ee
\end{enumerate}
\begin{enumerate}[$(\mathbf{H}_{A_2})$]
\item $0\in \text{Int}(\cD(A_2))$.
\end{enumerate}
\begin{enumerate}[$(\mathbf{H}^1_{b_{2}, \s_{2}})$]
\item There exists a constant $L_{b_{2}, \s_{2}}>0$ such that for $x_{i}\in\mR^n$, $\mu_{i}\in\cP_{2}(\mR^n)$, $i=1, 2$, $y\in\mR^m$,
\ce
&&|b_{2}(x_{1},\mu_1,y)-b_{2}(x_{2},\mu_2,y)|^{2}+\|\s_{2}(x_{1},\mu_1,y)-\s_{2}(x_{2},\mu_2,y)\|^{2}\\
&\leq& L_{b_{2}, \s_{2}}\(|x_{1}-x_{2}|^{2}+\mW^2_2(\mu_1,\mu_2)\).
\de
\end{enumerate}
\begin{enumerate}[$(\mathbf{H}^2_{b_{2}, \s_{2}})$]
\item For $(x,\mu)\in\mR^n\times\cP_{2}(\mR^n)$, $b_2(x,\mu,y)$ is continuous in $y$, and there exists a constant $\bar{L}_{b_{2}, \s_{2}}>0$ such that for $x\in\mR^n$, $\mu\in\cP_{2}(\mR^n)$, $y_i\in\mR^m, i=1,2$
\ce
&&2\<y_1-y_2, b_{2}(x,\mu,y_1)-b_{2}(x,\mu,y_2)\>\leq\bar{L}_{b_{2},\s_{2}}|y_{1}-y_{2}|^2,\\
&&\|\s_{2}(x,\mu,y_1)-\s_{2}(x,\mu,y_2)\|^2\leq\bar{L}_{b_{2},\s_{2}}|y_1-y_2|^2.
\de
\end{enumerate}
\begin{enumerate}[$(\mathbf{H}^3_{b_{2}, \s_{2}})$]
\item There exists a constant $\bar{\bar{L}}_{b_{2}, \s_{2}}>0$ such that for $x\in\mR^n$, $\mu\in\cP_{2}(\mR^n)$, $y\in\mR^m$, 
\be
|b_{2}(x,\mu,y)|^{2}+\|\s_{2}(x,\mu,y)\|^{2}
\leq \bar{\bar{L}}_{b_{2}, \s_{2}}(1+|x|^{2}+\mu(|\cdot|^2)+|y|^{2}).
\label{b2nu}
\ee
\end{enumerate}

\br
$(i)$ $(\mathbf{H}_{A_{2}})$ can be replaced by $\text{Int}(\cD(A_2))\ne\emptyset$. Indeed, if we require $\text{Int}(\cD(A_2))\ne\emptyset$, for any $a\in\text{Int}(\cD(A_2))$, by shifting the domain of $A_2$ and defining $\tilde{b}_2(x,\mu,y)= b_2(x,\mu,y-a), \tilde{\s}_2(x,\mu,y)= \s_2(x,\mu,y-a)$, this situation becomes the case of $0\in\text{Int}(\cD(A_2))$. However, we remind that $(\mathbf{H}_{A_{1}})$ can not be reduced to $0\in\text{Int}(\cD(A_1))$. That is because $b_1, \s_1$ depend on the distribution and the translation does not change the distribution. 

$(ii)$ Note that in $(\mathbf{H}^1_{b_{1}, \s_{1}})$ $(i)$, $b_1(x,\mu,y)$ is not Lipschitz continuous in $x$, and in $(\mathbf{H}^2_{b_{2}, \s_{2}})$, $b_2(x,\mu,y)$ is not Lipschitz continuous in $y$. Thus, $(\mathbf{H}^1_{b_{1}, \s_{1}})$ and $(\mathbf{H}^2_{b_{2}, \s_{2}})$ cover the Lipschitz continuous case.
\er

Now, it is the position to state the main result in this subsection.

\bt\label{well}
Assume that $(\mathbf{H}_{A_{1}})$, $(\mathbf{H}_{A_{2}})$, $(\mathbf{H}^1_{b_{1}, \s_{1}})$, $(\mathbf{H}^1_{b_{2}, \s_{2}})$-$(\mathbf{H}^3_{b_{2}, \s_{2}})$ hold. Then the system (\ref{Eqeu}) has a unique strong solution.
\et

The proof of the above theorem is placed in Section \ref{wellproo}.

\subsection{Averaging principles for multiscale multivalued McKean-Vlasov stochastic systems}

In this section, we present two averaging principles for multiscale multivalued McKean-Vlasov stochastic systems. 

Consider the system (\ref{Eqin1}), i.e.
\be\left\{\begin{array}{l}
\dif X_{t}^{\e,\d}\in -A_1(X_{t}^{\e,\d})\dif t+b_{1}(X_{t}^{\e,\d},\sL_{X_{t}^{\e,\d}},Y_{t}^{\e,\d})\dif t+\sqrt\e\s_{1}(X_{t}^{\e,\d},\sL_{X_{t}^{\e,\d}},Y_{t}^{\e,\d})\dif W^1_{t},\\
X_{0}^{\e,\d}=x_0\in\overline{\cD(A_1)},\quad  0\leq t\leq T,\\
\dif Y_{t}^{\e,\d}\in -A_2(Y_{t}^{\e,\d})\dif t+\frac{1}{\d}b_{2}(X_{t}^{\e,\d},\sL_{X_{t}^{\e,\d}},Y_{t}^{\e,\d})\dif t+\frac{1}{\sqrt{\d}}\s_{2}(X_{t}^{\e,\d},\sL_{X_{t}^{\e,\d}},Y_{t}^{\e,\d})\dif W^2_{t},\\
Y_{0}^{\e,\d}=y_0\in\overline{\cD(A_2)},\quad  0\leq t\leq T.
\end{array}
\right.
\label{Eqall}
\ee

Assume:
\begin{enumerate}[$(\mathbf{H}^{1'}_{b_{1}, \s_{1}})$]
\item
There exists a constant $L'_{b_{1}, \s_{1}}>0$ such that for $x_{i}\in\mR^n$, $\mu_{i}\in\cP_{2}(\mR^n)$, $y_{i}\in\mR^m$, $i=1, 2$,
\ce
&&|b_{1}(x_{1},\mu_1,y_{1})-b_{1}(x_{2},\mu_2,y_{2})|^{2}+\|\s_{1}(x_{1},\mu_1,y_1)-\s_{1}(x_{2},\mu_2,y_2)\|^{2}\\
&\leq& L'_{b_{1},\s_{1}}\(|x_{1}-x_{2}|^{2}+\mW^2_2(\mu_1,\mu_2)+|y_{1}-y_{2}|^{2}\).
\de
\end{enumerate}
\begin{enumerate}[$(\mathbf{H}^{2'}_{b_{2}, \s_{2}})$]
\item For $(x,\mu)\in\mR^n\times\cP_{2}(\mR^n)$, $b_2(x,\mu,y)$ is continuous in $y$, and there exist two constants $L'_{b_{2}, \s_{2}}, \b>0$ satisfying $\b>2L'_{b_2,\sigma_2}$ such that for $x\in\mR^n$, $\mu\in\cP_{2}(\mR^n)$, $y_{i}\in\mR^m$, $i=1, 2$,
\ce
\|\s_{2}(x,\mu,y_1)-\s_{2}(x,\mu,y_2)\|^2\leq L'_{b_{2},\s_{2}}|y_{1}-y_{2}|^2,
\de
\ce
2\<y_{1}-y_{2},b_{2}(x,\mu,y_{1})-b_{2}(x,\mu,y_{2})\>
+\|\s_{2}(x,\mu,y_{1})-\s_{2}(x,\mu,y_{2})\|^{2}\leq -\b|y_{1}-y_{2}|^{2}.
\de
\end{enumerate}

\br
$(i)$ $(\mathbf{H}^{1'}_{b_{1}, \s_{1}})$ is stronger than $(\mathbf{H}^{1}_{b_{1}, \s_{1}})$, and implies (\ref{b1line}). 
 
$(ii)$ $(\mathbf{H}^{2'}_{b_{2}, \s_{2}})$ is stronger than $(\mathbf{H}^2_{b_{2}, \s_{2}})$. Moreover, by $(\mathbf{H}^{2'}_{b_{2}, \s_{2}})$ $(\mathbf{H}^3_{b_{2}, \s_{2}})$, it holds that for $x\in\mR^n$, $\mu\in\cP_{2}(\mR^n)$, $y\in\mR^m$
\be
2\<y,b_{2}(x,\mu,y)\>+\|\s_{2}(x,\mu,y)\|^{2}\leq -\a|y|^{2}+C(1+|x|^{2}+\mu(|\cdot|^2)),
\label{bemu}
\ee
where $\a:=\b-2L'_{b_2,\sigma_2}>0$ and $C>0$ is a constant.
\er

By Theorem \ref{well}, under $(\mathbf{H}_{A_{1}})$, $(\mathbf{H}_{A_{2}})$, $(\mathbf{H}^{1'}_{b_{1}, \s_{1}})$, $(\mathbf{H}^1_{b_{2}, \s_{2}})$, $(\mathbf{H}^{2'}_{b_{2}, \s_{2}})$ and $(\mathbf{H}^3_{b_{2}, \s_{2}})$ we know that the system (\ref{Eqall}) has a unique strong solution still denoted as $(X_{\cdot}^{\e,\d},K_{\cdot}^{1,\e,\d},Y_{\cdot}^{\e,\d},K_{\cdot}^{2,\e,\d})$.

Take any $x\in \overline{\cD(A_1)}, \mu\in\cP_2(\overline{\cD(A_1)})$ and fix them. Consider the following multivalued SDE:
\be\left\{\begin{array}{l}
\dif Y_{t}^{x,\mu}\in -A_2(Y_{t}^{x,\mu})\dif t+b_{2}(x,\mu,Y_{t}^{x,\mu})\dif t+\s_{2}(x,\mu,Y_{t}^{x,\mu})\dif \hat{W}^2_{t},\\
Y_{0}^{x,\mu}=y_0\in\overline{\cD(A_2)}, \quad 0 \leq t \leq T,
\end{array}
\right.
\label{Eq2}
\ee
where $\hat{W}^2$ is a $d_2$-dimensional standard Brownian motion on another complete probability space $(\hat{\Omega}, \hat{\sF}, \hat{\mP})$. Under $(\mathbf{H}_{A_{2}})$ and $(\mathbf{H}^{2'}_{b_{2}, \s_{2}})$, we know that the above equation has a unique strong solution $(Y_{\cdot}^{x,\mu,y_0},K_{\cdot}^{2,x,\mu,y_0})$ (\cite{rwzx}). Moreover, by \cite[Theorem 3.2]{q2}, one could conclude that there exists a unique invariant probability measure $\nu^{x,\mu}$ for Eq.(\ref{Eq2}). 

Now, we state the first result in this subsection.
\bt \label{xbarxp}
Suppose that $(\mathbf{H}_{A_{1}})$, $(\mathbf{H}_{A_{2}})$, $(\mathbf{H}^{1'}_{b_{1}, \s_{1}})$, $(\mathbf{H}^{1}_{b_{2}, \s_{2}})$, $(\mathbf{H}^{2'}_{b_{2}, \s_{2}})$, $(\mathbf{H}^{3}_{b_{2}, \s_{2}})$ hold. If
\ce
\lim_{\e\rightarrow 0}\frac{\d}{\e}=\left\{\begin{array}{l}0,\\
\iota\in(0,\infty),
\end{array}
\right.
\de
it holds that
\ce
\lim_{\e\rightarrow 0}\mE\(\sup_{0\leq t\leq T}|X_{t}^{\e,\d}-\bar{X}^0_{t}|^{2}\)=0,
\de
where $(\bar{X}^{0},\bar{K}^{0})$ is a solution of the following equation:
\be\left\{\begin{array}{l}
\dif \bar{X}^0_{t}\in -A_1(\bar{X}^0_{t})\dif t+\bar{b}_{1}(\bar{X}^0_{t},D_{\bar{X}_t^0})\dif t,\\
\bar{X}^0_{0}=x_0\in\overline{\cD(A_1)},
\end{array}
\right.
\label{ldppsio0equde}
\ee
$\bar{b}_{1}(x,\mu):=\int_{\mR^{m}}b_{1}(x,\mu,y)\nu^{x,\mu}(\dif y)$ and $D_{\bar{X}_t^0}$ is the Dirac measure at $\bar{X}_t^0$.
\et

The proof of the above theorem is postponed to Section \ref{averprinproo}.

\br
If $A_1=0, A_2=0$, $b_1,b_2$ are independent of the distribution $\sL_{X_{t}^{\e,\d}}$, $\s_1, \s_2$ are $n, m$-order unit matrixes, respectively, $n=d_1, m=d_2$, and $\lim\limits_{\e\rightarrow 0}\frac{\d}{\e}=\iota\in(0,\infty)$, 
the system (\ref{Eqall}) is the same to the system $(1)+(2)$ with $s(\e)=1$ in \cite{abks}. There Athreya et al. proved that $\{X_{t}^{\e,\d}, 0\leq t\leq T\}$ converges in law on $C([0,T]; \mR^n)$ to $\{\bar{X}^0_{t}, 0\leq t\leq T\}$ (cf. \cite[Theorem 1.3]{abks}). Here in Theorem \ref{xbarxp}, we show that $\{X_{t}^{\e,\d}, 0\leq t\leq T\}$ converges in the mean square sense on $C([0,T]; \mR^n)$ to $\{\bar{X}^0_{t}, 0\leq t\leq T\}$. Therefore, our result is stronger.
\er

In the following, we take special $A_1$ and obtain a rate of $X^{\e,\d}$ converging to $\bar{X}^0$.

\bt\label{xbarxpcr}
Suppose that $A_1=\p I_{\cO}$, where $\cO$ is a closed and convex domain in $\mR^n$ with ${\rm Int}(\cO)\neq\emptyset$, and $(\mathbf{H}_{A_{2}})$, $(\mathbf{H}^{1'}_{b_{1}, \s_{1}})$, $(\mathbf{H}^{1}_{b_{2}, \s_{2}})$, $(\mathbf{H}^{2'}_{b_{2}, \s_{2}})$, $(\mathbf{H}^{3}_{b_{2}, \s_{2}})$ hold. If
\ce
\lim_{\e\rightarrow 0}\frac{\d}{\e}=\left\{\begin{array}{l}0,\\
\iota\in(0,\infty),
\end{array}
\right.
\de
it holds that for $0<\g<1$
\ce
\mE\(\sup_{0\leq t\leq T}|X_{t}^{\e,\d}-\bar{X}^0_{t}|^{2}\)\leq C(\d^{\g/2}+\d^\g+\d^{\frac{1}{2}(1-\g)}+\sqrt\e).
\de
\et

The proof of Theorem \ref{xbarxpcr} is placed in Section \ref{averprinproo}.

\br\label{conrat1}
If $\g=\frac{1}{3}$, by Theorem \ref{xbarxpcr}, it holds that
\ce
\mE\(\sup_{0\leq t\leq T}|X_{t}^{\e,\d}-\bar{X}^0_{t}|^{2}\)\leq C(\d^{\frac{1}{6}}+\sqrt\e).
\de
Note that in \cite[Theorem 2.1]{ghl}, for the following system 
\ce\left\{\begin{array}{l}
\dif X_{t}^{\e,\d}=b_{1}(X_{t}^{\e,\d},\sL_{X_{t}^{\e,\d}},Y_{t}^{\e,\d},\sL_{Y_{t}^{\e,\d}})\dif t+\sqrt{\e}\s_{1}(X_{t}^{\e,\d},\sL_{X_{t}^{\e,\d}},Y_{t}^{\e,\d},\sL_{Y_{t}^{\e,\d}})\dif W^1_{t},\\
X_{0}^{\e,\d}=x_0,\quad  0\leq t\leq T,\\
\dif Y_{t}^{\e,\d}=\frac{1}{\d}b_{2}(Y_{t}^{\e,\d})\dif t+\frac{1}{\sqrt{\d}}\s_{2}(Y_{t}^{\e,\d})\dif W^2_{t},\\
Y_{0}^{\e,\d}=y_0,\quad  0\leq t\leq T,
\end{array}
\right.
\de
Gao et al. obtained that
\ce
\mE\(\sup_{0\leq t\leq T}|X_{t}^{\e,\d}-\bar{X}^0_{t}|^{2}\)\leq C(\d^{\frac{1}{3}}+\e).
\de
That is, $A_1, A_2$ reduce the order of convergence.
\er

\subsection{The LDP for multiscale multivalued McKean-Vlasov stochastic systems}

In this subsection, we describe the LDP result for multiscale multivalued McKean-Vlasov stochastic systems. 

Here we require that $\s_1(x,\mu,y)=\s_1(x,\mu)$ and $(W^1_{t})$ and  $(W^2_{t})$ are mutually independent. Consider the following slow-fast system:
\be\left\{\begin{array}{l}
\dif X_{t}^{\e,\d}\in -A_1(X_{t}^{\e,\d})\dif t+b_{1}(X_{t}^{\e,\d},\sL_{X_{t}^{\e,\d}},Y_{t}^{\e,\d})\dif t+\sqrt{\e}\s_{1}(X_{t}^{\e,\d},\sL_{X_{t}^{\e,\d}})\dif W^1_{t},\\
X_{0}^{\e,\d}=x_0\in\overline{\cD(A_1)},\quad  0\leq t\leq T,\\
\dif Y_{t}^{\e,\d}\in -A_2(Y_{t}^{\e,\d})\dif t+\frac{1}{\d}b_{2}(X_{t}^{\e,\d},\sL_{X_{t}^{\e,\d}},Y_{t}^{\e,\d})\dif t+\frac{1}{\sqrt{\d}}\s_{2}(X_{t}^{\e,\d},\sL_{X_{t}^{\e,\d}},Y_{t}^{\e,\d})\dif W^2_{t},\\
Y_{0}^{\e,\d}=y_0\in\overline{\cD(A_2)},\quad  0\leq t\leq T,
\end{array}
\right.
\label{Eq1}
\ee
where $\d$ is a function of $\e$.

We assume more:
\begin{enumerate}[$(\mathbf{H}^4_{\s_{2}})$]
\item There exists a constant $L_{\sigma_2}>0$ such that for $x\in\mR^n$, $\mu\in\cP_{2}(\mR^n)$, 
\ce
\sup\limits_{y\in\mR^m}\|\s_2(x,\mu,y)\|^2\leq L_{\sigma_2}(1+|x|^2+\mu(|\cdot|^2)).
\de
\end{enumerate}

\br
$(\mathbf{H}^4_{\s_{2}})$ is weaker than the uniform boundedness.
\er

By Theorem \ref{well}, under $(\mathbf{H}_{A_{1}})$, $(\mathbf{H}_{A_{2}})$, $(\mathbf{H}^{1'}_{b_{1}, \s_{1}})$, $(\mathbf{H}^1_{b_{2}, \s_{2}})$, $(\mathbf{H}^{2'}_{b_{2}, \s_{2}})$ and $(\mathbf{H}^3_{b_{2}, \s_{2}})$ we know that the system (\ref{Eq1}) has a unique strong solution $(X_{\cdot}^{\e,\d},K_{\cdot}^{1,\e,\d},Y_{\cdot}^{\e,\d},K_{\cdot}^{2,\e,\d})$.

Now we present the main result in this subsection.

\bt\label{ldpmmsde}
Assume that $(\mathbf{H}_{A_{1}})$, $(\mathbf{H}_{A_{2}})$, $(\mathbf{H}^{1'}_{b_{1}, \s_{1}})$, $(\mathbf{H}^{1}_{b_{2}, \s_{2}})$, $(\mathbf{H}^{2'}_{b_{2}, \s_{2}})$, $(\mathbf{H}^{3}_{b_{2}, \s_{2}})$, $(\mathbf{H}^{4}_{\s_{2}})$ hold. If
\ce
\lim\limits_{\e\rightarrow 0}\frac{\d}{\e}=0,
\de
the family $\{X^{\e,\d},\e\in(0,1)\}$ satisfies the LDP in $\mS:=C([0,T],\overline{\mathcal{D}(A_1)})$ with the rate function given by
$$
I(\varsigma)=\frac{1}{2} \inf\limits_{h\in {\bf D}_{\varsigma}: \varsigma=\bar{X}^{h}}\|h\|_{\mH}^2,
$$
where $(\bar{X}^{h},\bar{K}^{h})$ solves the following equation
\ce\left\{\begin{array}{l}
\dif \bar{X}^h_{t}\in -A_1(\bar{X}^h_{t})\dif t+\bar{b}_{1}(\bar{X}^h_{t},D_{\bar{X}_t^0})\dif t+\s_{1}(\bar{X}^h_{t},D_{\bar{X}_t^0})\pi_1\dot{h}(t)\dif t,\\
\bar{X}^h_{0}=x_0\in\overline{\cD(A_1)},
\end{array}
\right.
\de
and $\pi_1: \mR^{d_1+d_2}\mapsto \mR^{d_1}$ is a projection operator.
\et

The proof of Theorem \ref{ldpmmsde} is placed in Section \ref{prooseco}.

\br
Based on the order that $\e, \d$ converge to $0$, there are three different regimes, i.e.
\ce
\lim_{\e\rightarrow 0}\frac{\d}{\e}=\left\{\begin{array}{l} 0, \qquad\qquad\qquad~~ ~\mbox{Regime}~ 1;\\
\iota\in (0,\infty), ~\mbox{Regime}~  2;\\
\infty, ~\qquad\quad\quad~\mbox{Regime}~ 3.
\end{array}
\right.
\de
Here we only can deal with Regime $1$. One reason lies in that we use the time discretization method. If $\lim\limits_{\e\rightarrow 0}\frac{\d}{\e}\neq 0$, some key estimates don't hold and we can not obtain the LDP for $\{X^{\e,\d},\e\in(0,1)\}$. The 
other reason is that if we follow up the method in \cite{ks2} to treat Regime $2$ and Regime $3$, it is difficult to express the needed Poisson equations.
\er

\section{Proof of Theorem \ref{well}}\label{wellproo}

In this section, we prove Theorem \ref{well} by a contractive mapping principle.

{\bf Proof of Theorem \ref{well}.}
First of all, let $\mS^2_{\mF}([0,T],\mR^m)$ be the space of all $\mF$-adapted processes $Y: \Omega\times [0,T]\mapsto\mR^m$ with $\mE[\sup\limits_{t\in[0,T]}|Y_t|^2]<\infty$. For any $\hat{Y}\in\mS^2_{\mF}([0,T],\mR^m)$, 
we consider the following system on $\mR^{n} \times \mR^{m}$:
\be\left\{\begin{array}{l}
\dif X_{t}\in -A_1(X_{t})\dif t+b_{1}(X_{t},\sL_{X_{t}},\hat{Y}_{t})\dif t+\s_{1}(X_{t},\sL_{X_{t}},\hat{Y}_{t})\dif W^1_{t},\\
X_{0}=\xi\in\overline{\cD(A_1)},\quad  0\leq t\leq T,\\
\dif Y_{t}\in -A_2(Y_{t})\dif t+b_{2}(X_{t},\sL_{X_{t}},Y_{t})\dif t+\s_{2}(X_{t},\sL_{X_{t}},Y_{t})\dif W^2_{t},\\
Y_{0}=\zeta\in\overline{\cD(A_2)},\quad  0\leq t\leq T.
\end{array}
\right.
\label{Eqeueu}
\ee
By \cite{gq}, we know that the $X$-equation in the system (\ref{Eqeueu}) has a unique strong solution $(X^{\hat{Y}}, K^{1,\hat{Y}})$ with $X^{\hat{Y}}\in\mS^2_{\mF}([0,T],\mR^n)$. Then inserting $X^{\hat{Y}}$ in the $Y$-equation, by \cite{rwzx} the $Y$-equation has a unique strong solution $(Y^{\hat{Y}}, K^{2,\hat{Y}})$ with $Y^{\hat{Y}}\in\mS^2_{\mF}([0,T],\mR^m)$. Define an operator $\sG: \mS^2_{\mF}([0,T],\mR^m)\mapsto\mS^2_{\mF}([0,T],\mR^m)$ by $\sG(\hat{Y})=Y^{\hat{Y}}$, and it is well defined.

Next, we show that $\sG$ is contractive. For any $\hat{Y}^1, \hat{Y}^2$, applying the It\^o formula to $|X^{1}_t-X^{2}_t|^2$, where $X^{1}_t:=X^{\hat{Y}^1}_t, X^{2}_t:=X^{\hat{Y}^2}_t$, we obtain that
\ce
|X^{1}_t-X^{2}_t|^2&=&-2\int_0^t\<X^{1}_s-X^{2}_s,\dif (K^{1,1}_s-K^{1,2}_s)\>\\
&&+2\int_0^t\<X^{1}_s-X^{2}_s, b_{1}(X^{1}_{s},\sL_{X^{1}_{s}},\hat{Y}^{1}_{s})-b_{1}(X^{2}_{s},\sL_{X^{2}_{s}},\hat{Y}^{2}_{s})\>\dif s\\
&&+2\int_0^t\<X^{1}_s-X^{2}_s,\(\s_{1}(X^{1}_{s},\sL_{X^{1}_{s}},\hat{Y}^{1}_{s})-\s_{1}(X^{2}_{s},\sL_{X^{2}_{s}},\hat{Y}^{2}_{s})\)\dif W_s^1\>\\
&&+\int_0^t\|\s_{1}(X^{1}_{s},\sL_{X^{1}_{s}},\hat{Y}^{1}_{s})-\s_{1}(X^{2}_{s},\sL_{X^{2}_{s}},\hat{Y}^{2}_{s})\|^2\dif s\\
&\leq&C\int_0^t\(|X^{1}_s-X^{2}_s|^2+\mW_2^2(\sL_{X^{1}_{s}},\sL_{X^{2}_{s}})\)\dif s\\
&&+\int_0^t|b_{1}(X^{2}_{s},\sL_{X^{2}_{s}},\hat{Y}^{1}_{s})-b_{1}(X^{2}_{s},\sL_{X^{2}_{s}},\hat{Y}^{2}_{s})|^2\dif s\\
&&+2\int_0^t\<X^{1}_s-X^{2}_s,\(\s_{1}(X^{1}_{s},\sL_{X^{1}_{s}},\hat{Y}^{1}_{s})-\s_{1}(X^{2}_{s},\sL_{X^{2}_{s}},\hat{Y}^{2}_{s})\)\dif W_s^1\>\\
&&+2\int_0^t\|\s_{1}(X^{2}_{s},\sL_{X^{2}_{s}},\hat{Y}^{1}_{s})-\s_{1}(X^{2}_{s},\sL_{X^{2}_{s}},\hat{Y}^{2}_{s})\|^2\dif s,
\de
where $K^{1,1}:=K^{\hat{Y}^1}, K^{1,2}:=K^{\hat{Y}^2}$. By taking the expectation on two sides and noticing the fact $\mW_2^2(\sL_{X^{1}_{s}},\sL_{X^{2}_{s}})\leq\mE|X^{1}_s-X^{2}_s|^2$, it holds that
\ce
\mE|X^{1}_t-X^{2}_t|^2\leq C\int_0^t\mE|X^{1}_s-X^{2}_s|^2\dif s+CT\mE\left[\sup\limits_{t\in[0,T]}|\hat{Y}^1_t-\hat{Y}^2_t|^2\right],
\de
which together with the Gronwall inequality yields that
\be
\mE|X^{1}_t-X^{2}_t|^2\leq CT\mE\left[\sup\limits_{t\in[0,T]}|\hat{Y}^1_t-\hat{Y}^2_t|^2\right]e^{CT}.
\label{xlxl1}
\ee

Besides, by the similar deduction to the above, it holds that
\be
\mE\sup\limits_{t\in[0,T]}|Y^{1}_t-Y^{2}_t|^2\leq C\int_0^T\mE\sup\limits_{r\in[0,s]}|Y^{1}_r-Y^{2}_r|^2\dif s+C\int_0^T\mE|X^{1}_s-X^{2}_s|^2\dif s,
\label{ylyl1}
\ee
where $Y^{1}_t:=Y^{\hat{Y}^1}_t, Y^{2}_t:=Y^{\hat{Y}^2}_t$. Inserting (\ref{xlxl1}) into (\ref{ylyl1}), by the Gronwall inequality we obtain that 
\ce
\mE\left[\sup\limits_{t\in[0,T]}|Y^{1}_t-Y^{2}_t|^2\right]\leq CT^2\mE\left[\sup\limits_{t\in[0,T]}|\hat{Y}^1_t-\hat{Y}^2_t|^2\right]e^{CT}.
\de
We take $0<T_0\leq T$ such that $CT_0^2e^{CT_0}<1$ and $\sG$ is contractive. Thus, there exists a unique fixed point $Y\in\mS^2_{\mF}([0,T_0],\mR^m)$ satisfying $\sG(Y)=Y$. Inserting $Y$ in the $X$-equation of the system (\ref{Eqeu}), we obtain 
that $(X, K^1, Y, K^2)$ is the unique solution for the system (\ref{Eqeu}). 

Finally, if $T_0=T$, the proof is complete; if $T_0<T$, we repeat the above procedure to get a unique solution of the system (\ref{Eqeu}) on $[T_0, T_1]$ for some $T_1\in(T_0, T]$. The approach is further utilized till $T_n=T$ so that a unique solution of the system (\ref{Eqeu}) on the whole interval $[0, T]$ is obtained. We are done.

\section{Proof of Theorem \ref{xbarxp} and \ref{xbarxpcr}}\label{averprinproo}

In this section, we prove Theorem \ref{xbarxp} and \ref{xbarxpcr}.

\subsection{Proof of Theorem \ref{xbarxp}}

In this subsection, we show Theorem \ref{xbarxp}. We begin with a lemma.

\bl
 Under $(\mathbf{H}_{A_{1}})$, $(\mathbf{H}_{A_{2}})$, $(\mathbf{H}^{1'}_{b_{1}, \s_{1}})$, $(\mathbf{H}^{1}_{b_{2}, \s_{2}})$, $(\mathbf{H}^{2'}_{b_{2}, \s_{2}})$, $(\mathbf{H}^{3}_{b_{2}, \s_{2}})$, there exists a constant $C>0$ such that 
\be
\mE\left(\sup\limits_{t\in[0,T]}|X_{t}^{\e,\d}|^{2}\right)\leq C(1+|x_0|^{2}+|y_0|^{2}).
\label{xeb}
\ee
\el
\begin{proof}
Note that $X_{t}^{\e,\d}$ satisfies the following equation:
\ce
X_{t}^{\e,\d}=x_0-K_t^{1,\e,\d}+\int_0^t b_{1}(X_{s}^{\e,\d},\sL_{X_{s}^{\e,\d}},Y_{s}^{\e,\d})\dif s+\sqrt \e\int_0^t\s_{1}(X_{s}^{\e,\d},\sL_{X_{s}^{\e,\d}},Y_{s}^{\e,\d})\dif W^1_{s}.
\de
For any $a\in\text{Int}(\cD(A_1))$, by applying the It\^o formula to $|X_{t}^{\e,\d}-a|^2$, it holds that 
\be
|X_{t}^{\e,\d}-a|^2&=&|x_0-a|^2-2\int_0^t\<X_{s}^{\e,\d}-a, \dif K_s^{1,\e,\d}\>+2\int_0^t\<X_{s}^{\e,\d}-a, b_{1}(X_{s}^{\e,\d},\sL_{X_{s}^{\e,\d}},Y_{s}^{\e,\d})\>\dif s\no\\
&&+2\sqrt\e\int_0^t\<X_{s}^{\e,\d}-a,\s_{1}(X_{s}^{\e,\d},\sL_{X_{s}^{\e,\d}},Y_{s}^{\e,\d})\dif W^1_{s}\>\no\\
&&+\e\int_0^t\|\s_{1}(X_{s}^{\e,\d},\sL_{X_{s}^{\e,\d}},Y_{s}^{\e,\d})\|^2\dif s, 
\label{xgito}
\ee
and furthermore 
\ce 
|X_{t}^{\e,\d}-a|^2&\leq&|x_0-a|^2-2M_1|K^{1,\e,\d}|_0^t+2M_2\int_0^t|X_{s}^{\e,\d}-a|\dif s+2M_3 t+\int_0^t|X_{s}^{\e,\d}-a|^2\dif s\\
&&+\int_0^t|b_{1}(X_{s}^{\e,\d},\sL_{X_{s}^{\e,\d}},Y_{s}^{\e,\d})|^2\dif s+2\sqrt\e\left|\int_0^t\<X_{s}^{\e,\d}-a,\s_{1}(X_{s}^{\e,\d},\sL_{X_{s}^{\e,\d}},Y_{s}^{\e,\d})\dif W^1_{s}\>\right|\\
&&+\e\int_0^t\|\s_{1}(X_{s}^{\e,\d},\sL_{X_{s}^{\e,\d}},Y_{s}^{\e,\d})\|^2\dif s,
\de
where Lemma \ref{inteineq} is used in the above inequality. Then the BDG inequality and $(\ref{b1line})$ imply that
\ce
\mE\left(\sup\limits_{s\in[0,t]}|X_{s}^{\e,\d}-a|^{2}\right)&\leq&(|x_0-a|^2+2(M_2+M_3)T)+(2M_2+1)\mE\int_0^t|X_{r}^{\e,\d}-a|^2\dif r\\
&&+\bar{L}_{b_{1}, \s_{1}}\mE\int_0^t(1+|X_{r}^{\e,\d}|^2+\sL_{X_{r}^{\e,\d}}(|\cdot|^2)+|Y_{r}^{\e,\d}|^2)\dif r\no\\
&&+2\mE\sup\limits_{s\in[0,t]}\left|\int_0^s\<X_{r}^{\e,\d}-a,\s_{1}(X_{r}^{\e,\d},\sL_{X_{r}^{\e,\d}},Y_{r}^{\e,\d})\dif W^1_{r}\>\right|\no\\ 
&\leq&C(|x_0-a|^2+1)+C\mE\int_{0}^{t}|X_{r}^{\e,\d}-a|^2\dif r+C\mE\int_{0}^{t}|Y_{r}^{\e,\d}|^{2}\dif r\no\\
&&+C\mE\left(\int_0^t|X_{r}^{\e,\d}-a|^2\|\s_{1}(X_{r}^{\e,\d},\sL_{X_{r}^{\e,\d}},Y_{r}^{\e,\d})\|^2\dif r\right)^{1/2}\no\\
&\leq&C(|x_0-a|^2+1)+C\mE\int_{0}^{t}|X_{r}^{\e,\d}-a|^2\dif r+C\mE\int_{0}^{t}|Y_{r}^{\e,\d}|^{2}\dif r\no\\
&&+\frac{1}{2}\mE\left(\sup\limits_{s\in[0,t]}|X_{s}^{\e,\d}-a|^{2}\right)+C\mE\int_0^t\|\s_{1}(X_{r}^{\e,\d},\sL_{X_{r}^{\e,\d}},Y_{r}^{\e,\d})\|^2\dif r,
\de
and furthermore
\be
\mE\left(\sup\limits_{s\in[0,t]}|X_{s}^{\e,\d}-a|^{2}\right)\leq C(|x_0-a|^2+1)+C\int_{0}^{t}\mE|X_{r}^{\e,\d}-a|^{2}\dif r+C\int_{0}^{t}\mE|Y_{r}^{\e,\d}|^{2}\dif r.
\label{exqc0}
\ee

For $Y_{t}^{\e,\d}$, applying the It\^{o} formula to $|Y_{t}^{\e,\d}|^{2}e^{\l t}$ for $\l=\frac{\a}{2\d}$ and taking the expectation, one could obtain that for any $v\in A_2(0)$
\ce
\mE|Y_{t}^{\e,\d}|^{2}e^{\l t}&=& |y_0|^{2}+\l\mE\int_0^t|Y_{s}^{\e,\d}|^{2}e^{\l s}\dif s-2\mE\int_0^te^{\l s}\<Y_{s}^{\e,\d},\dif K_s^{2,\e,\d}\>\\
&&+\frac{2}{\d}\mE\int_{0}^{t}e^{\l s}\<Y_{s}^{\e,\d}, b_{2}(X_{s}^{\e,\d},\sL_{X_{s}^{\e,\d}},Y_{s}^{\e,\d})\>\dif s\\
&&+\frac{1}{\d}\mE\int_{0}^{t}e^{\l s}\|\s_{2}(X_{s}^{\e,\d},\sL_{X_{s}^{\e,\d}},Y_{s}^{\e,\d})\|^2\dif s\\
&\overset{(\ref{bemu})}{\leq}& |y_0|^{2}+\l\mE\int_0^t|Y_{s}^{\e,\d}|^{2}e^{\l s}\dif s+2\mE\int_0^te^{\l s}|v||Y_{s}^{\e,\d}|\dif s\\
&&+\frac{1}{\d}\mE\int_{0}^{t}e^{\l s}\(-\a|Y_{s}^{\e,\d}|^{2}+C(1+|X_{s}^{\e,\d}|^{2}+\sL_{X_{s}^{\e,\d}}(|\cdot|^2))\)\dif s\\
&\leq& |y_0|^{2}+\(\l+\frac{\a}{2\d}-\frac{\a}{\d}\)\mE\int_0^t|Y_{s}^{\e,\d}|^{2}e^{\l s}\dif s+\frac{2\d}{\a}|v|^2\int_0^te^{\l s}\dif s\\
&&+\frac{C}{\d}\mE\int_{0}^{t}e^{\l s}(1+|X_{s}^{\e,\d}|^{2}+\mE|X_{s}^{\e,\d}|^{2})\dif s\\
&\leq& |y_0|^{2}+\frac{2|v|^2}{\a}\frac{e^{\l t}-1}{\l}+C\frac{e^{\l t}-1}{\d\l}\left(\mE\left(\sup\limits_{s\in[0,t]}|X_{s}^{\e,\d}|^{2}\right)+1\right),
\de
where we use Lemma \ref{equi} in the first inequality and in the second inequality the fact that
\ce
2|v||Y_{s}^{\e,\d}|\leq \frac{2\d}{\a}|v|^2+\frac{\a}{2\d}|Y_{s}^{\e,\d}|^{2}.
\de
From this, it follows that
\be
\mE|Y_{t}^{\e,\d}|^{2}\leq C(|y_0|^{2}+1)+C\mE\left(\sup\limits_{s\in[0,t]}|X_{s}^{\e,\d}|^{2}\right).
\label{yges}
\ee

Inserting (\ref{yges}) in (\ref{exqc0}), by the Gronwall inequality one can get that
\ce
\mE\left(\sup\limits_{t\in[0,T]}|X_{t}^{\e,\d}|^{2}\right)&\leq& 2\mE\left(\sup\limits_{t\in[0,T]}|X_{t}^{\e,\d}-a|^{2}\right)+2|a|^2\\
&\leq& C(1+|x_0-a|^2+|y_0|^{2})+2|a|^2\\
&\leq& C(1+|x_0|^2+|y_0|^{2}).
\de
The proof is complete.
\end{proof}

Next, by the similar deduction to that for \cite[Lemma 5.3]{q1}, we have the following result. 

\bl\label{xehatxets}
Under $(\mathbf{H}_{A_{1}})$, $(\mathbf{H}_{A_{2}})$, $(\mathbf{H}^{1'}_{b_{1}, \s_{1}})$, $(\mathbf{H}^{1}_{b_{2}, \s_{2}})$, $(\mathbf{H}^{2'}_{b_{2}, \s_{2}})$, $(\mathbf{H}^{3}_{b_{2}, \s_{2}})$, we have that
\be
\lim\limits_{l\rightarrow 0}\sup\limits_{s\in[0,T]}\mE\sup _{s \leqslant t \leqslant s+l}|X_{t}^{\e,\d}-X_{s}^{\e,\d}|^{2}=0.
 \label{xegts}
\ee
\el

Besides, we recall Eq.(\ref{ldppsio0equde}), i.e.
\ce\left\{\begin{array}{l}
\dif \bar{X}^0_{t}\in -A_1(\bar{X}^0_{t})\dif t+\bar{b}_{1}(\bar{X}^0_{t},D_{\bar{X}_t^0})\dif t,\\
\bar{X}^0_{0}=x_0\in\overline{\cD(A_1)}.
\end{array}
\right.
\de
The following lemma guarantees the well-posedness of Eq.(\ref{ldppsio0equde}).

\bl \label{averc}
Under $(\mathbf{H}_{A_{1}})$, $(\mathbf{H}_{A_{2}})$, $(\mathbf{H}^{1'}_{b_{1}, \s_{1}})$, $(\mathbf{H}^{1}_{b_{2}, \s_{2}})$, $(\mathbf{H}^{2'}_{b_{2}, \s_{2}})$, $(\mathbf{H}^{3}_{b_{2}, \s_{2}})$, Eq.(\ref{ldppsio0equde}) has a unique solution $(\bar{X}^0_{\cdot},\bar{K}^0_{\cdot})$. Moreover, it holds that
\be
&&\sup\limits_{t\in[0,T]}|\bar{X}^0_{t}|^{2}\leq C(1+|x_0|^{2}), \label{barx0b}\\
&&\lim\limits_{l\rightarrow 0}\sup\limits_{s\in[0,T]}\sup\limits_{s\leq t\leq s+l}|\bar{X}^0_{t}-\bar{X}^0_{s}|^2=0. \label{barx0ts}
\ee
\el
\begin{proof}
First of all, by the similar deduction to that of \cite[Lemma 3.8]{rsx}, we know that for any $x_1, x_2\in\overline{\cD(A_1)}$ and $\mu_1,\mu_2\in\cP_2(\overline{\cD(A_1)})$
\be
|\bar{b}_1(x_1,\mu_1)-\bar{b}_1(x_2,\mu_2)|\leq C(|x_1-x_2|+\mW_2(\mu_1,\mu_2)),
\label{barb1lip}
\ee
which implies that Eq.(\ref{ldppsio0equde}) has a unique solution $(\bar{X}^0_{\cdot},\bar{K}^0_{\cdot})$. Then, since the proofs of the required estimates are similar to that for $X^{\e,\d}$, we omit them. The proof is complete.
\end{proof}

{\bf Proof of Theorem \ref{xbarxp}.} Since the proof of the above theorem is similar to that for \cite[Theorem 3.1]{q1}, we only give a sketch. First of all, we construct two auxiliary processes $\hat{X}^{\e,\d}, \hat{Y}^{\e,\d}$. Then by the above results we estimate 
$\mE\sup\limits_{t\in[0,T]}|X_{t}^{\e,\d}-\hat{X}_t^{\e,\d}|^2$ and $\mE\sup\limits_{t\in[0,T]}|\hat{X}_t^{\e,\d}-\bar{X}^0_{t}|^2$, respectively. Note that 
$$
\mE\sup\limits_{t\in[0,T]}|X^{\e,\d}_t-\bar{X}^0_{t}|^2\leq 2\mE\sup\limits_{t\in[0,T]}|X_{t}^{\e,\d}-\hat{X}_t^{\e,\d}|^2+2\mE\sup\limits_{t\in[0,T]}|\hat{X}_t^{\e,\d}-\bar{X}^0_{t}|^2.
$$
Thus, we obtain the required result.

\subsection{Proof of Theorem \ref{xbarxpcr}}\label{prootheo5.6}

In this subsection, we prove Theorem \ref{xbarxpcr}. We prepare a key lemma.

\bl\label{xehatxetscr}
Under the assumptions of Theorem \ref{xbarxpcr}, we have that for $l>0$ small enough 
\be
&&\sup\limits_{s\in[0,T]}\mE\sup _{s \leqslant t \leqslant s+l}|X_{t}^{\e,\d}-X_{s}^{\e,\d}|^{2}\leq Cl, \label{xegtscr}\\
&&\sup\limits_{s\in[0,T]}\sup\limits_{s\leq t\leq s+l}|\bar{X}^0_{t}-\bar{X}^0_{s}|^2\leq Cl, \label{barx0tscr}
\ee
where the constant $C>0$ is independent of $\e, \d, l$.
\el
\begin{proof}
Since the proofs of (\ref{xegtscr}) and (\ref{barx0tscr}) are similar, we only prove (\ref{xegtscr}).

First of all, note that for $0\leq s\leq t\leq s+l\leq T$, 
\ce
X_{t}^{\e,\d}-X_{s}^{\e,\d}=-K_t^{1,\e,\d}+K_s^{1,\e,\d}+\int_{s}^{t}b_{1}(X_{r}^{\e,\d},\sL_{X_r^{\e,\d}},Y_{r}^{\e,\d})\dif r+\sqrt\e\int_{s}^{t}\s_{1}(X_{r}^{\e,\d},\sL_{X_r^{\e,\d}},Y_{r}^{\e,\d})\dif W^1_{r}.
\de
Thus, the It\^o formula, (\ref{b1line}) and \cite[Lemma 2.4]{arrst} imply that
\ce
|X_{t}^{\e,\d}-X_{s}^{\e,\d}|^2&=&-2\int_s^t\<X_{r}^{\e,\d}-X_{s}^{\e,\d},\dif K_r^{1,\e,\d}\>\no\\
&&+2\int_s^t\<X_{r}^{\e,\d}-X_{s}^{\e,\d},b_{1}(X_{r}^{\e,\d},\sL_{X_r^{\e,\d}},Y_{r}^{\e,\d})\>\dif r\no\\
&&+2\sqrt\e\int_s^t\<X_{r}^{\e,\d}-X_{s}^{\e,\d},\s_{1}(X_{r}^{\e,\d},\sL_{X_r^{\e,\d}},Y_{r}^{\e,\d})\dif W^1_{r}\>\no\\
&&+\e\int_s^t\|\s_{1}(X_{r}^{\e,\d},\sL_{X_r^{\e,\d}},Y_{r}^{\e,\d})\|^2\dif r\no\\
&\leq&\int_s^t|X_{r}^{\e,\d}-X_{s}^{\e,\d}|^2\dif r+C\int_s^t\(1+|X_{r}^{\e,\d}|^2+\sL_{X_r^{\e,\d}}(|\cdot|^2)+|Y_{r}^{\e,\d}|^2\)\dif r\\
&&+2\left|\int_s^t\<X_{r}^{\e,\d}-X_{s}^{\e,\d},\s_{1}(X_{r}^{\e,\d},\sL_{X_r^{\e,\d}},Y_{r}^{\e,\d})\dif W^1_{r}\>\right|.
\de
By the BDG inequality, it holds that
\ce
\mE\sup\limits_{t\in[s,s+l]}|X_{t}^{\e,\d}-X_{s}^{\e,\d}|^2&\leq& \int_s^{s+l}\mE|X_{r}^{\e,\d}-X_{s}^{\e,\d}|^2\dif r+C\int_s^{s+l}(1+\mE|X_{r}^{\e,\d}|^2+\mE|Y_{r}^{\e,\d}|^2)\dif r\\
&&+C\mE\(\int_s^{s+l}|X_{r}^{\e,\d}-X_{s}^{\e,\d}|^2\|\s_{1}(X_{r}^{\e,\d},\sL_{X_r^{\e,\d}},Y_{r}^{\e,\d})\|^2\dif r\)^{1/2}\\
&\leq& \int_s^{s+l}\mE\sup\limits_{u\in[s,r]}|X_{u}^{\e,\d}-X_{s}^{\e,\d}|^2\dif r+Cl+\frac{1}{2}\mE\sup\limits_{t\in[s,s+l]}|X_{t}^{\e,\d}-X_{s}^{\e,\d}|^2\\
&&+C\mE\int_s^{s+l}\|\s_{1}(X_{r}^{\e,\d},\sL_{X_r^{\e,\d}},Y_{r}^{\e,\d})\|^2\dif r\\
&\leq& \int_s^{s+l}\mE\sup\limits_{u\in[s,r]}|X_{u}^{\e,\d}-X_{s}^{\e,\d}|^2\dif r+Cl+\frac{1}{2}\mE\sup\limits_{t\in[s,s+l]}|X_{t}^{\e,\d}-X_{s}^{\e,\d}|^2,
\de
which together with the Gronwall inequality yields (\ref{xegtscr}). The proof is complete.
\end{proof} 

{\bf Proof of Theorem \ref{xbarxpcr}.} We replace (\ref{xegts}), (\ref{barx0ts}) by (\ref{xegtscr}), (\ref{barx0tscr}), respectively, and show Theorem \ref{xbarxpcr} by the same deduction to that of Theorem \ref{xbarxp}.

\section{Proof of Theorem \ref{ldpmmsde}}\label{prooseco}

In this section, we prove Theorem \ref{ldpmmsde} by a weak convergence approach. We  divide the proof into three parts. In the first subsection, we construct the mapping $\cG^\e$ and make some estimates. In the second subsection, we define the mapping $\cG^0$ and give a few estimates. Verification for Condition \ref{cond} is placed in the third subsection. 

\subsection{Construction of $\cG^\e$ and related estimates}\label{some}

In this subsection, we construct the mapping $\cG^\e$ and make some estimates.

By Theorem \ref{well}, under the assumptions of Theorem \ref{ldpmmsde}, we know that the system (\ref{Eq1}) has a unique strong solution $(X_{\cdot}^{\e,\d},K_{\cdot}^{1,\e,\d},Y_{\cdot}^{\e,\d},K_{\cdot}^{2,\e,\d})$. Thus, there exists a functional $\cG^\e=\cG_{\sL_{X^{\e,\d}}}^{\e}: C([0,T];\mathbb{R}^{d_1+d_2})\mapsto C([0,T],\overline{\cD(A_1)})$ such that 
\ce
X^{\e,\d}=\cG^{\e}(\sqrt{\e}W), \quad W:=(W^1,W^2).
\de

Next, we consider the following controlled processes:
\be\left\{\begin{array}{l}
\dif X_{t}^{\e,\d,u}\in -A_1(X_{t}^{\e,\d,u})\dif t+b_{1}(X_{t}^{\e,\d,u},\sL_{X_t^{\e,\d}},Y_{t}^{\e,\d,u})\dif t+\sigma_1(X^{\e,\d,u}_{t},\sL_{X_t^{\e,\d}})\pi_1\dot{u}(t)\dif t\\
\qquad\qquad+\sqrt{\e}\s_{1}(X_{t}^{\e,\d,u},\sL_{X_t^{\e,\d}})\dif W^1_{t},\\
X_{0}^{\e,\d,u}=x_0\in\overline{\cD(A_1)},\quad  0\leq t\leq T,\\
\dif Y_{t}^{\e,\d,u}\in -A_2(Y_{t}^{\e,\d,u})\dif t+\frac{1}{\d}b_{2}(X_{t}^{\e,\d,u},\sL_{X_t^{\e,\d}},Y_{t}^{\e,\d,u})\dif t\\
\qquad\qquad +\frac{1}{\sqrt{\d \e}}\sigma_2(X^{\e,\d,u}_{t},\sL_{X_t^{\e,\d}},Y^{\e,\d,u}_{t})\pi_2\dot{u}(t)\dif t+\frac{\sqrt{\e}}{\sqrt{\d \e}}\s_{2}(X_{t}^{\e,\d,u},\sL_{X_t^{\e,\d}},Y_{t}^{\e,\d,u})\dif W^2_{t},  \\
Y_{0}^{\e,\d,u}=y_0\in\overline{\cD(A_2)},\quad  0\leq t\leq T,\quad u\in\mathbf{A}_2^{N},
\end{array}
\right.
\label{contproc}
\ee
where $\pi_1: \mR^{d_1+d_2}\mapsto \mR^{d_1}, \pi_2: \mR^{d_1+d_2}\mapsto \mR^{d_2}$ are two projection operators. Thus, by the Girsanov theorem, the system $(\ref{contproc})$ have a unique strong solution denoted by $(X^{\e,\d,u}, K^{1,\e,\d,u},Y^{\e,\d,u}, K^{2,\e,\d,u})$. Moreover, $X^{\e,\d,u}=\cG^{\e}(\sqrt{\e}W+u)$. 

\br\label{chmv}
Here we emphasize that the controlled system $(\ref{contproc})$ contains the distribution of $X_t^{\e,\d}$ but not the distribution of $X_{t}^{\e,\d,u}$, which is a key characteristic of McKean-Vlasov stochastic systems.
\er

\bl \label{xutzutc}
Under the assumptions of Theorem \ref{ldpmmsde}, for $\{u_{\e}, \e\in(0,1)\}\subset\mathbf{A}_{2}^{N}$, there exists a constant $C>0$ such that 
\be
&&\mE\left(\sup\limits_{t\in[0,T]}|X_{t}^{\e,\d,u_\e}|^{2}\right)\leq C(1+|x_0|^{2}+|y_0|^{2}), \label{xeub}\\
&&\int_0^T\mE|Y_{r}^{\e,\d,u_\e}|^{2}\dif r\leq C(1+|x_0|^{2}+|y_0|^{2}),\label{yeub}\\
&&\mE|K^{1,\e,\d,u_\e}|_0^T\leq C(1+|x_0|^{2}+|y_0|^{2}). \label{keub}
\ee
\el
\begin{proof}
First of all, we estimate $X_{t}^{\e,\d,u_\e}$. Note that $X_{t}^{\e,\d,u_\e}$ satisfies the following equation:
\ce
X_{t}^{\e,\d,u_\e}&=&x_0-K_t^{1,\e,\d,u_\e}+\int_0^t b_{1}(X_{s}^{\e,\d,u_\e},\sL_{X_s^{\e,\d}},Y_{s}^{\e,\d,u_\e})\dif s+\int_{0}^{t}\sigma_1(X^{\e,\d,u_\e}_{s},\sL_{X_s^{\e,\d}})\pi_1\dot{u}_\e(s)\dif s\\
&&+\sqrt{\e}\int_0^t\s_{1}(X_{s}^{\e,\d,u_\e},\sL_{X_s^{\e,\d}})\dif W^1_{s}.
\de
The It\^o formula yields that for any $a\in\text{Int}(\cD(A_1))$
\be
&&|X_{t}^{\e,\d,u_\e}-a|^2\no\\
&=&|x_0-a|^2-2\int_0^t\<X_{s}^{\e,\d,u_\e}-a,\dif K_s^{1,\e,\d,u_\e}\>+2\int_0^t\<X_{s}^{\e,\d,u_\e}-a, b_{1}(X_{s}^{\e,\d,u_\e},\sL_{X_s^{\e,\d}},Y_{s}^{\e,\d,u_\e})\>\dif s\no\\
&&+2\int_0^t\<X_{s}^{\e,\d,u_\e}-a,\sigma_1(X^{\e,\d,u_\e}_{s},\sL_{X_s^{\e,\d}})\pi_1\dot{u}_\e(s)\>\dif s+\e\int_0^t\|\s_{1}(X_{s}^{\e,\d,u_\e},\sL_{X_s^{\e,\d}})\|^2\dif s\no\\
&&+2\sqrt{\e}\int_0^t\<X_{s}^{\e,\d,u_\e}-a,\s_{1}(X_{s}^{\e,\d,u_\e},\sL_{X_s^{\e,\d}})\dif W^1_{s}\>.
\label{itoxegu}
\ee
Then the BDG inequality, Lemma \ref{inteineq} and $(\ref{b1line})$ imply that 
\ce
&&\mE\left(\sup\limits_{s\in[0,t]}|X_{s}^{\e,\d,u_\e}-a|^{2}\right)\\
&\leq&C(|x_0-a|^{2}+1)+C\mE\int_{0}^{t}|X_{r}^{\e,\d,u_\e}-a|^2\dif r+C\int_{0}^{t}\sL_{X_r^{\e,\d}}(|\cdot|^2)\dif r+C\mE\int_{0}^{t}|Y_{r}^{\e,\d,u_\e}|^{2}\dif r\\
&&+\frac{1}{4}\mE\left(\sup\limits_{s\in[0,t]}|X_{s}^{\e,\d,u_\e}-a|^{2}\right)+C\mE\left(\int_0^t\|\sigma_1(X^{\e,\d,u_\e}_{r},\sL_{X_r^{\e,\d}})\|^2\dif r\right)\left(\int_0^t|\dot{u}_\e(r)|^2\dif r\right)\no\\
&&+C\mE\left(\int_0^t|X_{r}^{\e,\d,u_\e}-a|^2\|\s_{1}(X_{r}^{\e,\d,u_\e},\sL_{X_r^{\e,\d}})\|^2\dif r\right)^{1/2}\no\\
&\leq&C(|x_0-a|^{2}+1)+C\mE\int_{0}^{t}|X_{r}^{\e,\d,u_\e}-a|^2\dif r+C\int_{0}^{t}\mE|X_r^{\e,\d}|^2\dif r+C\mE\int_{0}^{t}|Y_{r}^{\e,\d,u_\e}|^{2}\dif r\no\\
&&+\frac{1}{2}\mE\left(\sup\limits_{s\in[0,t]}|X_{s}^{\e,\d,u_\e}-a|^{2}\right),
\de
and furthermore by (\ref{xeb})
\be
\mE\left(\sup\limits_{s\in[0,t]}|X_{s}^{\e,\d,u_\e}-a|^{2}\right)&\leq& C(1+|x_0|^{2}+|y_0|^{2})+C\int_{0}^{t}\mE|X_{r}^{\e,\d,u_\e}-a|^{2}\dif r\no\\
&&+C\int_{0}^{t}\mE|Y_{r}^{\e,\d,u_\e}|^{2}\dif r.
\label{exqcu}
\ee

For $Y_{t}^{\e,\d,u_\e}$, fix $v\in A_2(0)$. Applying the It\^{o} formula to $|Y_{t}^{\e,\d,u_\e}|^{2}e^{\l t}$ for $\l=\frac{\a}{3\d}$ and taking the expectation, one could obtain that
\ce
\mE|Y_{t}^{\e,\d,u_\e}|^{2}e^{\l t}
&=& |y_0|^{2}+\l\mE\int_0^t|Y_{s}^{\e,\d,u_\e}|^{2}e^{\l s}\dif s-2\mE\int_0^te^{\l s}\<Y_{s}^{\e,\d,u_\e},\dif K_s^{2,\e,\d,u_\e}\>\\
&&+\frac{2}{\d}\mE\int_{0}^{t}e^{\l s}\<Y_{s}^{\e,\d,u_\e}, b_{2}(X_{s}^{\e,\d,u_\e},\sL_{X_s^{\e,\d}},Y_{s}^{\e,\d,u_\e})\>\dif s\\
&&+\frac{1}{\d}\mE\int_{0}^{t}e^{\l s}\|\s_{2}(X_{s}^{\e,\d,u_\e},\sL_{X_s^{\e,\d}},Y_{s}^{\e,\d,u_\e})\|^2\dif s\\
&&+\frac{2}{\sqrt{\d \e}}\mE\int_{0}^{t}e^{\l s}\<Y_{s}^{\e,\d,u_\e},\sigma_2(X^{\e,\d,u_\e}_{s},\sL_{X_s^{\e,\d}},Y^{\e,\d,u_\e}_{s})\pi_2\dot{u}_\e(s)\>\dif s\\
&\overset{(\ref{bemu})}{\leq}& |y_0|^{2}+\l\mE\int_0^t|Y_{s}^{\e,\d,u_\e}|^{2}e^{\l s}\dif s+2\mE\int_0^te^{\l s}|v||Y_{s}^{\e,\d,u_\e}|\dif s\\
&&+\frac{1}{\d}\mE\int_{0}^{t}e^{\l s}\left(-\a|Y_{s}^{\e,\d,u_\e}|^{2}+C\(1+|X_{s}^{\e,\d,u_\e}|^{2}+\sL_{X_s^{\e,\d}}(|\cdot|^2)\)\right)\dif s\\
&&+\frac{2}{\sqrt{\d \e}}\mE\int_{0}^{t}e^{\l s}|Y_{s}^{\e,\d,u_\e}|\|\sigma_2(X^{\e,\d,u_\e}_{s},\sL_{X_s^{\e,\d}},Y^{\e,\d,u_\e}_{s})\||\dot{u}_\e(s)|\dif s.
\de
Note that for the third and last terms of the right side for the above inequality
\ce
2|v||Y_{s}^{\e,\d,u_\e}|\leq \frac{\a}{3\d}|Y_{s}^{\e,\d,u_\e}|^2+\frac{3\d}{\a}|v|^2,
\de
and
\ce
&&\frac{2}{\sqrt{\d \e}}|Y_{s}^{\e,\d,u_\e}|\|\sigma_2(X^{\e,\d,u_\e}_{s},\sL_{X_s^{\e,\d}},Y^{\e,\d,u_\e}_{s})\||\dot{u}_\e(s)|\\
&\leq&\frac{\a}{3\d}|Y_{s}^{\e,\d,u_\e}|^{2}+\frac{C}{\e}\(1+|X^{\e,\d,u_\e}_{s}|^2+\sL_{X_s^{\e,\d}}(|\cdot|^2)\)|\dot{u}_\e(s)|^2,
\de
where $(\mathbf{H}^4_{\s_{2}})$ is used. Thus, we have that
\ce
\mE|Y_{t}^{\e,\d,u_\e}|^{2}e^{\l t}&\leq& |y_0|^{2}+C\int_0^te^{\l s}\dif s+\(\l+\frac{2\a}{3\d}-\frac{\a}{\d}\)\mE\int_0^t|Y_{s}^{\e,\d,u_\e}|^{2}e^{\l s}\dif s\\
&&+\frac{C}{\d}\mE\int_{0}^{t}e^{\l s}\(1+|X_{s}^{\e,\d,u_\e}|^{2}+\sL_{X_s^{\e,\d}}(|\cdot|^2)\)\dif s\\
&&+\frac{C}{\e}\mE\int_{0}^{t}e^{\l s}\(1+|X^{\e,\d,u_\e}_{s}|^2+\sL_{X_s^{\e,\d}}(|\cdot|^2)\)|\dot{u}_\e(s)|^2\dif s\\
&\leq& |y_0|^{2}+C\frac{e^{\l t}-1}{\l}+C\frac{e^{\l t}-1}{\d\l}\left(1+\mE\sup\limits_{s\in[0,t]}|X_{s}^{\e,\d,u_\e}|^{2}+\mE\sup\limits_{s\in[0,t]}|X_{s}^{\e,\d}|^{2}\right)\\
&&+\frac{C}{\e}\mE\(1+\sup\limits_{s\in[0,T]}|X^{\e,\d,u_\e}_{s}|^2+\sup\limits_{s\in[0,T]}\sL_{X_s^{\e,\d}}(|\cdot|^2)\)\int_{0}^{t}e^{\l s}|\dot{u}_\e(s)|^2\dif s.
\de
From this, it follows that
\ce
\mE|Y_{t}^{\e,\d,u_\e}|^{2}&\leq& C(1+|x_0|^{2}+|y_0|^{2})+C\mE\left(\sup\limits_{s\in[0,t]}|X_{s}^{\e,\d,u_\e}|^{2}\right)\\
&&+\frac{C}{\e}\mE\(1+\sup\limits_{s\in[0,T]}|X^{\e,\d,u_\e}_{s}|^2+\sup\limits_{s\in[0,T]}\sL_{X_s^{\e,\d}}(|\cdot|^2)\)\int_{0}^{t}e^{-\l (t-s)}|\dot{u}_\e(s)|^2\dif s,
\de
and furthermore
\be
\int_0^t\mE|Y_{r}^{\e,\d,u_\e}|^{2}\dif r&\leq& CT(1+|x_0|^{2}+|y_0|^{2})+C\int_0^t\mE\left(\sup\limits_{s\in[0,r]}|X_{s}^{\e,\d,u_\e}|^{2}\right)\dif r\no\\
&&+\frac{C}{\e}\mE\(1+\sup\limits_{s\in[0,T]}|X^{\e,\d,u_\e}_{s}|^2+\sup\limits_{s\in[0,T]}\sL_{X_s^{\e,\d}}(|\cdot|^2)\)\int_0^t\int_{0}^{r}e^{-\l (r-s)}|\dot{u}_\e(s)|^2\dif s\dif r\no\\
&\leq& CT(1+|x_0|^{2}+|y_0|^{2})+C\int_0^t\mE\left(\sup\limits_{s\in[0,r]}|X_{s}^{\e,\d,u_\e}|^{2}\right)\dif r\no\\
&&+C\left(\frac{\d}{\e}\right)\mE\(1+\sup\limits_{s\in[0,T]}|X^{\e,\d,u_\e}_{s}|^2+\sup\limits_{s\in[0,T]}\sL_{X_s^{\e,\d}}(|\cdot|^2)\)\int_0^t|\dot{u}_\e(s)|^2\dif s\no\\
&\leq& C(1+|x_0|^{2}+|y_0|^{2})+C\int_0^t\mE\left(\sup\limits_{s\in[0,r]}|X_{s}^{\e,\d,u_\e}|^{2}\right)\dif r\no\\
&&+CN\frac{\d}{\e}\mE\sup\limits_{s\in[0,T]}|X^{\e,\d,u_\e}_{s}|^2.
\label{zeues}
\ee
where we use $u_\e\in {\bf A}_2^N$ and $\lim\limits_{\e\rightarrow0}\frac{\d}{\e}=0$ in the last inequality. 

Inserting (\ref{zeues}) in (\ref{exqcu}), one can get that
\ce
\mE\left(\sup\limits_{s\in[0,t]}|X_{s}^{\e,\d,u_\e}-a|^{2}\right)&\leq& C(1+|x_0|^{2}+|y_0|^{2})+C\int_{0}^{t}\mE\left(\sup\limits_{s\in[0,r]}|X_{s}^{\e,\d,u_\e}-a|^{2}\right)\dif r\\
&&+CN\frac{\d}{\e}\mE\sup\limits_{s\in[0,T]}|X^{\e,\d,u_\e}_{s}-a|^2.
\de
Since $\lim\limits_{\e\rightarrow0}\frac{\d}{\e}=0$, there exists a $\e_0>0$ such that for $\e<\e_0$, $CN\frac{\d}{\e}<1$ and 
$$
\mE\left(\sup\limits_{s\in[0,t]}|X_{s}^{\e,\d,u_\e}-a|^{2}\right)\leq C(1+|x_0|^{2}+|y_0|^{2})+C\int_{0}^{t}\mE\left(\sup\limits_{s\in[0,r]}|X_{s}^{\e,\d,u_\e}-a|^{2}\right)\dif r,
$$
which together with the Gronwall inequality yields (\ref{xeub}). Then by (\ref{xeub}) and (\ref{zeues}), we obtain (\ref{yeub}).

Finally, for $K^{1,\e,\d,u_\e}$, by (\ref{itoxegu}) and Lemma \ref{inteineq}, it holds that
\ce
&&2M_1\mE\left| K^{1,\e,\d,u_\e} \right|_{0}^{T}\\
&\leq& |x_0-a|^2+2(M_2+M_3)T+(2M _2+1)\int_0^T\mE{\left| X^{\e,\d,u_\e}_s-a\right|^2}\dif s+\int_0^T\mE|X_{s}^{\e,\d,u_\e}|^2\dif s\\
&&+\mE\left(\sup\limits_{s\in[0,T]}|X_{s}^{\e,\d,u_\e}-a|^2\right)+C\int_0^T\(1+\mE{\left| X^{\e,\d,u_\e}_s\right|^2}+\sL_{X_s^{\e,\d}}(|\cdot|^2)+\mE{\left| Y^{\e,\d,u_\e}_s\right|^2}\)\dif s,
\de
which together with (\ref{xeub}) and (\ref{yeub}) yields the required estimate. The proof is complete.
\end{proof}

Moreover, by the similar deduction to that for \cite[Lemma 5.3]{q1}, we obtain the following limit result.

\bl
Under the assumptions of Theorem \ref{ldpmmsde}, for $\{u_{\e}, \e\in(0,1)\}\subset\mathbf{A}_{2}^{N}$, we have that
\be
\lim\limits_{l\rightarrow 0}\sup\limits_{s\in[0,T]}\mE\sup _{s \leqslant t \leqslant s+l}|X_{t}^{\e,\d,u_\e}-X_{s}^{\e,\d,u_\e}|^{2}=0.
\label{xegutse}
\ee
\el

Finally we introduce the following auxiliary process: for any $u_{\e}\in\mathbf{A}_{2}^{N}$,
\be\left\{\begin{array}{l}
\hat{Y}_{t}^{\e,\d,u_\e}\in -A_2(\hat{Y}_{t}^{\e,\d,u_\e})\dif t+\frac{1}{\d}b_{2}(X_{k\triangle}^{\e,\d,u_\e},\sL_{X_{k\triangle}^{\e,\d}},\hat{Y}_{t}^{\e,\d,u_\e})\dif t\\
\qquad\qquad +\frac{1}{\sqrt{\d}}\s_{2}(X_{k\triangle}^{\e,\d,u_\e},\sL_{X_{k\triangle}^{\e,\d}},\hat{Y}_{t}^{\e,\d,u_\e})\dif W^2_{t}, \quad t\in[k\triangle,(k+1)\triangle),\\
\hat{Y}_{k\triangle}^{\e,\d,u_\e}=Y_{k\triangle}^{\e,\d,u_\e}, \quad \hat{K}_{k\triangle}^{2,\e,\d,u_\e}=K_{k\triangle}^{2,\e,\d,u_\e},
\end{array}
\right.
\label{hatzu}
\ee
where $\triangle$ is a fixed positive number depending on $\d$, and $k=0,1,\cdots,[\frac{T}{\triangle}]$, and $[\frac{T}{\triangle}]$ denotes the integer part of $\frac{T}{\triangle}$. Moreover, we mention the fact that $[\frac{t}{\triangle}]=k$ for $t\in[k\triangle, (k+1)\triangle)$.

\bl
Under the assumptions of Theorem \ref{ldpmmsde}, for $\{u_{\e}, \e\in(0,1)\}\subset\mathbf{A}_{2}^{N}$, there exists a constant $C>0$ such that 
 \be
 \sup\limits_{t\in[0,T]}\mE|\hat{Y}_{t}^{\e,\d,u_\e}|^{2}&\leq& C(1+|x_0|^{2}+|y_0|^{2}), \label{hatzub}\\
\int_0^T\mE|Y_{t}^{\e,\d,u_\e}-\hat{Y}_{t}^{\e,\d,u_\e}|^{2}\dif t&\leq& CN\left(\frac{\d}{\e}\right)+\frac{CT}{\a}\left(\sup\limits_{s\in[0,T]}\mE\sup _{s \leqslant r \leqslant s+\triangle}|X_{r}^{\e,\d,u_\e}-X_{s}^{\e,\d,u_\e}|^{2}\right)\no\\
&&+\frac{CT}{\a}\left(\sup\limits_{s\in[0,T]}\mE\sup _{s \leqslant r \leqslant s+\triangle}|X_{r}^{\e,\d}-X_{s}^{\e,\d}|^{2}\right).
\label{unztu}
\ee
\el
\begin{proof}
Since the proof of (\ref{hatzub}) is standard, we only prove (\ref{unztu}).

First of all, by (\ref{contproc}) and (\ref{hatzu}), we have that for $t\in[k\triangle,(k+1)\triangle)$
\ce
&&Y_{t}^{\e,\d,u_\e}-\hat{Y}_{t}^{\e,\d,u_\e}\\
&=&-K_t^{2,\e,\d,u_\e}+\hat{K}_t^{2,\e,\d,u_\e}+\frac{1}{\d}\int_{k\triangle}^{t}\(b_{2}(X_{s}^{\e,\d,u_\e},\sL_{X_{s}^{\e,\d}},Y_{s}^{\e,\d,u_\e})
-b_{2}(X_{k\triangle}^{\e,\d,u_\e},\sL_{X_{k\triangle}^{\e,\d}},\hat{Y}_{s}^{\e,\d,u_\e})\)\dif s\\
&&+\frac{1}{\sqrt{\d}}\int_{k\triangle}^{t}\(\s_{2}(X_{s}^{\e,\d,u_\e},\sL_{X_{s}^{\e,\d}},Y_{s}^{\e,\d,u_\e})
-\s_{2}(X_{k\triangle}^{\e,\d,u_\e},\sL_{X_{k\triangle}^{\e,\d}},\hat{Y}_{s}^{\e,\d,u_\e})\)\dif W^2_{s}\\
&&+\frac{1}{\sqrt{\d \e}}\int_{k\triangle}^{t}\s_{2}(X_{s}^{\e,\d,u_\e},\sL_{X_{s}^{\e,\d}},Y_{s}^{\e,\d,u_\e})\pi_2\dot{u}_\e(s)\dif s.
\de
Applying the It\^{o} formula to $|Y_{t}^{\e,\d,u_\e}-\hat{Y}_{t}^{\e,\d,u_\e}|^{2}e^{\l t}$ for $\l=\frac{\a}{2\d}$ and taking the expectation, by $(\mathbf{H}^{1}_{b_{2},\s_{2}})$ and $(\mathbf{H}^{2'}_{b_{2},\s_{2}})$ one could obtain that
\ce
&&\mE|Y_{t}^{\e,\d,u_\e}-\hat{Y}_{t}^{\e,\d,u_\e}|^{2}e^{\l t}\\
&\leq&\l\mE\int_{k\triangle}^{t}|Y_{s}^{\e,\d,u_\e}-\hat{Y}_{s}^{\e,\d,u_\e}|^{2}e^{\l s}\dif s\\
&&+\frac{1}{\d}\mE\int_{k \triangle}^{t}2e^{\l s}\<Y_{s}^{\e,\d,u_\e}-\hat{Y}_{s}^{\e,\d,u_\e}, b_{2}(X_{s}^{\e,\d,u_\e},\sL_{X_{s}^{\e,\d}},Y_{s}^{\e,\d,u_\e})-b_{2}(X_{s}^{\e,\d,u_\e},\sL_{X_{s}^{\e,\d}},\hat{Y}_{s}^{\e,\d,u_\e})\>\dif s\\
&&+\frac{1}{\d}\mE\int_{k \triangle}^{t}e^{\l s}\|\s_{2}(X_{s}^{\e,\d,u_\e},\sL_{X_{s}^{\e,\d}},Y_{s}^{\e,\d,u_\e})-\s_{2}(X_{s}^{\e,\d,u_\e},\sL_{X_{s}^{\e,\d}},\hat{Y}_{s}^{\e,\d,u_\e})\|^{2}\dif s\\
&&+\frac{\a}{2\d}\mE\int_{k \triangle}^{t}|Y_{s}^{\e,\d,u_\e}-\hat{Y}_{s}^{\e,\d,u_\e}|^{2}e^{\l s}\dif s+\frac{C}{\e}\mE\int_{k\triangle}^{t}e^{\l s}\(1+|X_{s}^{\e,\d,u_\e}|^2+\sL_{X_{s}^{\e,\d}}(|\cdot|^2)\)|\dot{u}_\e(s)|^2\dif s\\
&&+\frac{1}{\d}\mE\int_{k \triangle}^{t}2e^{\l s}\<Y_{s}^{\e,\d,u_\e}-\hat{Y}_{s}^{\e,\d,u_\e}, b_{2}(X_{s}^{\e,\d,u_\e},\sL_{X_{s}^{\e,\d}},\hat{Y}_{s}^{\e,\d,u_\e})-b_{2}(X_{k \triangle}^{\e,\d,u_\e},\sL_{X_{k \triangle}^{\e,\d}},\hat{Y}_{s}^{\e,\d,u_\e})\>\dif s\\
&&+\frac{1}{\d}\mE\int_{k\triangle}^{t}e^{\l s}\|\s_{2}(X_{s}^{\e,\d,u_\e},\sL_{X_{s}^{\e,\d}},Y_{s}^{\e,\d,u_\e})-\s_{2}(X_{s}^{\e,\d,u_\e},\sL_{X_{s}^{\e,\d}},\hat{Y}_{s}^{\e,\d,u_\e})\|^{2}\dif s\\
&&+\frac{1}{\d}\mE\int_{k\triangle}^{t}2e^{\l s}\|\s_{2}(X_{s}^{\e,\d,u_\e},\sL_{X_{s}^{\e,\d}},\hat{Y}_{s}^{\e,\d,u_\e})-\s_{2}(X_{k \triangle}^{\e,\d,u_\e},\sL_{X_{k \triangle}^{\e,\d}},\hat{Y}_{s}^{\e,\d,u_\e})\|^{2}\dif s\\
&\leq&(\l-\frac{\b}{\d}+\frac{2L'_{b_2,\s_2}}{\d}+\frac{\a}{2\d})\mE\int_{k\triangle}^{t}|Y_{s}^{\e,\d,u_\e}-\hat{Y}_{s}^{\e,\d,u_\e}|^{2}e^{\l s}\dif s\\
&&+\frac{C}{\e}\mE\int_{k\triangle}^{t}e^{\l s}\(1+|X_{s}^{\e,\d,u_\e}|^2+\sL_{X_{s}^{\e,\d}}(|\cdot|^2)\)|\dot{u}_\e(s)|^2\dif s\\
&&+\frac{C}{\d}\mE\int_{k\triangle}^{t}e^{\l s}|X_{s}^{\e,\d,u_\e}-X_{k \triangle}^{\e,\d,u_\e}|^{2}\dif s+\frac{C}{\d}\int_{k\triangle}^{t}e^{\l s}\mW_2^2(\sL_{X_{s}^{\e,\d}},\sL_{X_{k \triangle}^{\e,\d}})\dif s\\
&\leq&\frac{C}{\e}\mE\left(1+\sup\limits_{s\in[0,T]}|X_{s}^{\e,\d,u_\e}|^2+\sup\limits_{s\in[0,T]}\sL_{X_{s}^{\e,\d}}(|\cdot|^2)\right)\int_{k\triangle}^{t}e^{\l s}|\dot{u}_\e(s)|^2\dif s\\
&&+\frac{C}{\d}\left(\sup\limits_{s\in[0,T]}\mE\sup _{s \leqslant r \leqslant s+\triangle}|X_{r}^{\e,\d,u_\e}-X_{s}^{\e,\d,u_\e}|^{2}\right)\frac{e^{\l t}-e^{\l k\triangle}}{\l}\\
 &&+\frac{C}{\d}\left(\sup\limits_{s\in[0,T]}\mE\sup _{s \leqslant r \leqslant s+\triangle}|X_{r}^{\e,\d}-X_{s}^{\e,\d}|^{2}\right)\frac{e^{\l t}-e^{\l k \triangle}}{\l}.
\de

Finally, it follows that
\ce
\int_{k\triangle}^{t}\mE|Y_{r}^{\e,\d,u_\e}-\hat{Y}_{r}^{\e,\d,u_\e}|^{2}\dif r&\leq& C\left(\frac{\d}{\e}\right)\mE\left(1+\sup\limits_{s\in[0,T]}|X_{s}^{\e,\d,u_\e}|^2+\sup\limits_{s\in[0,T]}\sL_{X_{s}^{\e,\d}}(|\cdot|^2)\right)\int_{k\triangle}^{t}|\dot{u}_\e(s)|^2\dif s\\
&&+\frac{C}{\a}\left(\sup\limits_{s\in[0,T]}\mE\sup _{s \leqslant r \leqslant s+\triangle}|X_{r}^{\e,\d,u_\e}-X_{s}^{\e,\d,u_\e}|^{2}\right)\triangle\\
&&+\frac{C}{\a}\left(\sup\limits_{s\in[0,T]}\mE\sup _{s \leqslant r \leqslant s+\triangle}|X_{r}^{\e,\d}-X_{s}^{\e,\d}|^{2}\right)\triangle,
\de
which implies the required estimate. The proof is complete.
\end{proof}

\subsection{Construction of $\cG^0$ and related estimates}

In this subsection, we define the mapping $\cG^0$ and give a few estimates.

Consider the following multivalued differential equation:
\be\left\{\begin{array}{l}
\dif\bar{X}^{u}_{t}\in A_1(\bar{X}^{u}_{t})\dif t+\bar{b}_{1}(\bar{X}^u_{t},D_{\bar{X}^0_{t}})\dif t+\s_{1}(\bar{X}^u_{t},D_{\bar{X}^0_{t}})\pi_1\dot{u}(t)\dif t, \quad u\in\mathbf{A}_2^{N},\\
\bar{X}^{u}_{0}=x_0\in\overline{\cD(A_1)}.
\end{array}
\right.
\label{deteequa}
\ee
By (\ref{barb1lip}) and $(\mathbf{H}^{1'}_{b_{1}, \s_{1}})$, it holds that Eq.(\ref{deteequa}) has a unique solution $(\bar{X}^{u}, \bar{K}^{u})$. Define a map $\cG^{0}: \mH\mapsto C([0,T],\overline{\cD(A_1)})$ by $\cG^{0}(u)=\bar{X}^{u}$, and we make some key estimates.

\bl
Under the assumptions of Theorem \ref{ldpmmsde}, it holds that
\be
&&\sup\limits_{t\in[0,T]}|\bar{X}^u_{t}|^{2}\leq C(1+|x_0|^{2}), a.s.\label{barxub}\\
&&|\bar{K}^{u}|_0^T\leq C(1+|x_0|^{2}), a.s. \label{barkub}\\
&&\lim\limits_{l\rightarrow 0}\sup\limits_{s\in[0,T]}\sup\limits_{s\leq t\leq s+l}|\bar{X}^u_{t}-\bar{X}^u_{s}|^2=0, a.s.. \label{barxuts}
\ee
\el

Since the proofs of (\ref{barxub}), (\ref{barkub}) and (\ref{barxuts}) are similar to that of (\ref{xeub}), (\ref{keub}) and (\ref{xegutse}), respectively, we omit them.

\subsection{Verification for Condition \ref{cond}}

In this subsection, we verify Condition \ref{cond} $(i)$ and $(ii)$.

\bl\label{auxilemm2}
Suppose that the assumptions of Theorem \ref{ldpmmsde} hold, and $h_{\e}\rightarrow h$ in $\mathbf{D}_2^{N}$ as ${\e}\rightarrow0$. Then $\cG^{0}(h_\e)$ converges to $\cG^{0}(h)$.
\el
\begin{proof}
By the definition of $\cG^{0}$, $\cG^{0}(h_{\e})$ and $\cG^{0}(h)$ satisfy the following equations respectively:
\ce
&&\bar{X}^{h_\e}_{t}=x_0-\bar{K}^{h_\e}_{t}+\int_{0}^{t}\bar{b}_1(\bar{X}^{h_\e}_{s}, D_{\bar{X}^{0}_{s}})\dif s+\int_{0}^{t}\sigma_1(\bar{X}^{h_\e}_{s},D_{\bar{X}^{0}_{s}})\pi_1\dot{h}_{{\e}}(s)\dif s,\\
&&\bar{X}^{h}_{t}=x_0-\bar{K}^{h}_{t}+\int_{0}^{t}\bar{b}_1(\bar{X}^{h}_{s}, D_{\bar{X}^{0}_{s}})\dif s+\int_{0}^{t}\sigma_1(\bar{X}^{h}_{s},D_{\bar{X}^{0}_{s}})\pi_1\dot{h}(s)\dif s.
\de
Set $Z^{0}(t)=\bar{X}^{h_{\e}}_{t}-\bar{X}^{h}_{t}$, and by Lemma \ref{equi} and (\ref{barb1lip}) we have
\be
|Z^{0}(t)|^{2}&=&-2\int_{0}^{t}\langle Z^{0}(s),\dif(\bar{K}^{h_{\e}}_{s}-\bar{K}^{h}_{s})\rangle+2\int_{0}^{t}\langle Z^{0}(s),\bar{b}_1(\bar{X}^{h_{\e}}_{s},D_{\bar{X}^{0}_{s}})-\bar{b}_1(\bar{X}^{h}_{s},D_{\bar{X}^{0}_{s}})\rangle\dif s  \no\\
&&+2\int_{0}^{t}\langle Z^{0}(s),\sigma_1(\bar{X}^{h_{\e}}_{s},D_{\bar{X}^0_{s}})\pi_1\dot{h}_{{\e}}(s)-\sigma_1(\bar{X}^{h}_{s},D_{\bar{X}^0_{s}})\pi_1\dot{h}(s)\rangle\dif s  \no\\
&\leq&C\int_{0}^{t}|Z^{0}(s)|^2\dif s+2\int_{0}^{t}\langle Z^{0}(s),\(\sigma_1(\bar{X}^{h_{\e}}_{s},D_{\bar{X}^0_{s}})- \sigma_1(\bar{X}^{h}_{s},D_{\bar{X}^0_{s}})\)\pi_1\dot{h}_{{\e}}(s)\rangle \dif s\no\\
&&+2\int_{0}^{t}\langle Z^{0}(s),\sigma_1(\bar{X}^{h}_{s},D_{\bar{X}^0_{s}})(\pi_1\dot{h}_{{\e}}(s)-\pi_1\dot{h}(s))\rangle \dif s\no\\
&=:&C\int_{0}^{t}|Z^{0}(s)|^2\dif s+I_1(t)+I_2(t).
\label{eq10}
\ee

Next, for $I_1(t)$, it follows from $(\mathbf{H}^{1'}_{b_{1}, \s_{1}})$  and $h_\e\in\mathbf{D}_2^{N}$ that
\be
\sup\limits_{s\in[0,t]}|I_1(s)|&\leq&2\sup\limits_{s\in[0,t]}\left|\int_{0}^{s}\langle Z^{0}(r),\(\sigma_1(\bar{X}^{h_\e}_{r},D_{\bar{X}^0_{r}})- \sigma_1(\bar{X}^{h}_{r},D_{\bar{X}^0_{r}})\)\pi_1\dot{h}_{{\e}}(r)\rangle \dif r\right|\no\\
&\leq& 2\sqrt{L_{b_1,\s_1}}\left(\int_{0}^{t}|Z^{0}(r)|^{4}\dif r\right)^{\frac{1}{2}}\left(\int_{0}^{t}|\dot{h}_{{\e}}(r)|^{2}\dif r\right)^{\frac{1}{2}}\no\\
&\leq& \frac{1}{2}\sup\limits_{r\in[0,t]}|Z^{0}(r)|^{2}+C\int_{0}^{t}|Z^{0}(r)|^{2}\dif r.
\label{eq11}
\ee

Inserting (\ref{eq11}) into (\ref{eq10}), we obtain that
\ce
\sup\limits_{s\in[0,t]}|Z^{0}(t)|^{2}\leq \frac{1}{2}\sup\limits_{s\in[0,t]}|Z^{0}(s)|^{2}+C\int_{0}^{t}|Z^{0}(s)|^{2}\dif s+\sup\limits_{s\in[0,T]}|I_2(s)|,
\de
which together with Gronwall's inequality implies that
\be
\sup\limits_{t\in[0,T]}|\bar{X}^{h_{{\e}}}_{t}-\bar{X}^{h}_{t}|^2\leq \sup\limits_{s\in[0,T]}|I_2(s)|e^{CT}.
\label{xhexh}
\ee

In the following, we are devoted to proving that $\lim\limits_{{\e}\rightarrow0}\sup\limits_{s\in[0,T]}|I_2(s)|=0$. First, set for $h_\e, h\in\mathbf{D}_{2}^{N}$
\ce
g_\e(t):=\int_0^t\sigma_1(\bar{X}^{h}_{r},D_{\bar{X}^0_{r}})(\pi_1\dot{h}_{{\e}}(r)-\pi_1\dot{h}(r))\dif r,
\de
and it holds that 
\be
\lim\limits_{\e\rightarrow0}\sup\limits_{t\in[0,T]}|g_\e(t)|=0.
\label{gees}
\ee
The proof of (\ref{gees}) is postponed in the Appendix.

For $I_{2}(t)$, applying the integral by parts formula to $\<Z^0(s),g_\e(s)\>$, we have that
\ce
\frac{1}{2}I_{2}(t)&=&\<Z^0(t),g_\e(t)\>+\int_0^t\<g_\e(s),\dif (\bar{K}^{h_{\e}}_{s}-\bar{K}^{h}_{s})\>\\
&&-\int_0^t\<g_\e(s),\bar{b}_1(\bar{X}^{h_{\e}}_{s},D_{\bar{X}^{0}_{s}})-\bar{b}_1(\bar{X}^{h}_{s},D_{\bar{X}^{0}_{s}})\>\dif s\\
&&-\int_0^t\<g_\e(s),\sigma_1(\bar{X}^{h_{\e}}_{s},D_{\bar{X}^0_{s}})\pi_1\dot{h}_{{\e}}(s)-\sigma_1(\bar{X}^{h}_{s},D_{\bar{X}^0_{s}})\pi_1\dot{h}(s)\>\dif s\\
&=:&I_{21}(t)+I_{22}(t)+I_{23}(t)+I_{24}(t).
\de
For $I_{21}(t)$, note that 
\ce
\sup\limits_{t\in[0,T]}|I_{21}(t)|\leq\sup\limits_{t\in[0,T]}|Z^0(t)|\sup\limits_{t\in[0,T]}|g_\e(t)|.
\de
Thus, by (\ref{barxub}) and (\ref{gees}), it holds that
\ce
\sup\limits_{t\in[0,T]}|I_{21}(t)|\rightarrow 0, \quad \e\rightarrow 0.
\de
For $I_{22}(t)$, noticing that
\ce
\sup\limits_{t\in[0,T]}|I_{22}(t)|\leq\sup\limits_{t\in[0,T]}|g_\e(t)|(|\bar{K}^{h_\e}|_0^T+|\bar{K}^{h}|_0^T),
\de
by (\ref{barkub}) and (\ref{gees}) we obtain that $\sup\limits_{t\in[0,T]}|I_{22}(t)|\rightarrow 0$. By the same deduction to the above, one can get that
\ce
\lim\limits_{\e\rightarrow 0}\sup\limits_{t\in[0,T]}|I_{23}(t)|=0, \quad \lim\limits_{\e\rightarrow 0}\sup\limits_{t\in[0,T]}|I_{24}(t)|=0.
\de

Combining the above deduction, we get that 
\be
\lim\limits_{\e\rightarrow 0}\sup\limits_{t\in[0,T]}\left|I_{2}(t)\right|=0.
\label{j22}
\ee

By (\ref{xhexh}) and (\ref{j22}), it holds that 
\ce
\lim\limits_{\e\rightarrow 0}\sup\limits_{t\in[0,T]}\left|\cG^{0}\left(h_{\e}\right)(t)-\cG^{0}\left(h\right)(t)\right|=\lim\limits_{\e\rightarrow 0}\sup\limits_{t\in[0,T]}|\bar{X}^{h_{\e}}_{t}-\bar{X}^{h}_{t}|^2=0,
\de
which completes the proof. 
\end{proof}

\bl\label{auxilemm3}
Suppose that the assumptions of Theorem \ref{ldpmmsde} hold. Assume that $\{u_{\e},\epsilon>0\}\subset \mathbf{A}_{2}^{N}$. Then for any $\eta>0$,
$$
\lim\limits_{\e\rightarrow0}\mP\left(\sup\limits_{t\in[0,T]}\left|\cG^{\e}\left(\sqrt{\e}W_\cdot+u_{\e}\right)(t)-\cG^{0}\left(u_\e\right)(t)\right|>\eta\right)=0.
$$
\el
\begin{proof}
We divide the proof into two steps. In the first step, we estimate 
$$
\sup\limits_{t\in[0,T]}\left|\cG^{\e}\left(\sqrt{\e}W_\cdot+u_{\e}\right)(t)-\cG^{0}\left(u_\e\right)(t)\right|. 
$$
In the second step, we show the required result.

{\bf Step 1.} We estimate 
$$
\sup\limits_{t\in[0,T]}\left|\cG^{\e}\left(\sqrt{\e}W_\cdot+u_{\e}\right)(t)-\cG^{0}\left(u_\e\right)(t)\right|. 
$$

Note that
$$
X^{\e,\d,u_{\e}}=\cG^{\e}\left(\sqrt{\e}W+u_{\e}\right), \quad \bar{X}^{u_\e}=\cG^{0}\left(u_\e\right).
$$
Thus, set $Z^{\e,u_{\e}}(t)=X^{\e,\d,u_{\e}}_{t}-\bar{X}^{u_\e}_{t}$, and it holds that
\ce
Z^{\e,u_{\e}}(t)
&=&-(K^{1,\e,\d,u_{\e}}_{t}-\bar{K}^{u_\e}_{t})+\int_{0}^{t}\left[b_1(X^{\e,\d,u_{\e}}_{s},\sL_{X_s^{\e,\d}},Y^{\e,\d,u_{\e}}_{s})-\bar{b}_1(\bar{X}^{u_\e}_{s},D_{\bar{X}^{0}_{s}})\right]\dif s\no\\
&&+\int_{0}^{t}\left[\s_1(X^{\e,\d,u_{\e}}_{s},\sL_{X_s^{\e,\d}})\pi_1\dot{u}_{\e}(s)-\sigma_1(\bar{X}^{u_\e}_{s},D_{\bar{X}^{0}_{s}})\pi_1\dot{u}_\e(s)\right]\dif s\\
&&+\sqrt{\e}\int_{0}^{t}\s_1(X^{\e,\d,u_{\e}}_{s},\sL_{X_s^{\e,\d}})\dif W^1_s.
\de
By It\^o's formula and Lemma \ref{equi}, we get that
\be
|Z^{\e,u_{\e}}(t)|^{2}
&\leq&2\int_{0}^{t}\langle   Z^{\e,u_{\e}}(s), b_1(X^{\e,\d,u_{\e}}_{s},\sL_{X_s^{\e,\d}},Y^{\e,\d,u_{\e}}_{s})-\bar{b}_1(\bar{X}^{u_\e}_{s},D_{\bar{X}^{0}_{s}}) \rangle  \dif s   \no\\
&&+2\int_{0}^{t}\langle   Z^{\e,u_{\e}}(s), \s_1(X^{\e,\d,u_{\e}}_{s},\sL_{X_s^{\e,\d}})\pi_1\dot{u}_{\e}(s)-\sigma_1(\bar{X}^{u_\e}_{s},D_{\bar{X}^{0}_{s}})\pi_1\dot{u}_\e(s) \rangle  \dif s  \no\\
&&+2\sqrt{\e} \int_{0}^{t}\langle  Z^{\e,u_{\e}}(s),  \s_1(X^{\e,\d,u_{\e}}_{s},\sL_{X_s^{\e,\d}}) \dif W^1_s\rangle +\e\int_{0}^{t}\|\s_1(X^{\e,\d,u_{\e}}_{s},\sL_{X_s^{\e,\d}})\|^{2} \dif s\no\\
&=:&J_{1}(t)+J_{2}(t)+J_{3}(t)+J_{4}(t).
\label{j1j2j3j4}
\ee

For $J_{1}(t)$, note that
\ce
J_{1}(t)&=&2\int_{0}^{t}\langle   Z^{\e,u_{\e}}(s), b_1(X^{\e,\d,u_{\e}}_{s},\sL_{X_s^{\e,\d}},Y^{\e,\d,u_{\e}}_{s})-b_1(X^{\e,\d,u_{\e}}_{s(\triangle)},\sL_{X_{s(\triangle)}^{\e,\d}},\hat{Y}^{\e,\d,u_{\e}}_{s}) \rangle  \dif s\\
&&+2\int_{0}^{t}\langle   Z^{\e,u_{\e}}(s),-\bar{b}_1(X^{\e,\d,u_{\e}}_{s},\sL_{X_s^{\e,\d}})+\bar{b}_1(X^{\e,\d,u_{\e}}_{s(\triangle)},\sL_{X_{s(\triangle)}^{\e,\d}})\rangle \dif s\\
&&+2\int_{0}^{t}\langle   Z^{\e,u_{\e}}(s),\bar{b}_1(X^{\e,\d,u_{\e}}_{s},\sL_{X_{s}^{\e,\d}})-\bar{b}_1(\bar{X}^{u_\e}_{s},D_{\bar{X}^{0}_{s}}) \rangle  \dif s\\
&&+2\int_{0}^{t}\langle   Z^{\e,u_{\e}}(s)-Z^{\e,u_{\e}}(s(\triangle)),b_1(X^{\e,\d,u_{\e}}_{s(\triangle)},\sL_{X_{s(\triangle)}^{\e,\d}},\hat{Y}^{\e,\d,u_{\e}}_{s})-\bar{b}_1(X^{\e,\d,u_{\e}}_{s(\triangle)},\sL_{X_{s(\triangle)}^{\e,\d}})\rangle \dif s\\
&&+2\int_{0}^{t}\langle   Z^{\e,u_{\e}}(s(\triangle)),b_1(X^{\e,\d,u_{\e}}_{s(\triangle)},\sL_{X_{s(\triangle)}^{\e,\d}},\hat{Y}^{\e,\d,u_\e}_{s})-\bar{b}_1(X^{\e,\d,u_{\e}}_{s(\triangle)},\sL_{X_{s(\triangle)}^{\e,\d}})\rangle \dif s\\
&=:&J_{11}(t)+J_{12}(t)+J_{13}(t)+J_{14}(t)+J_{15}(t),
\de
where $s(\triangle):=[\frac{s}{\triangle}]\triangle$. So, by the H\"older inequality and the Lipschitz continuity of $b_1, \bar{b}_1$, we get that
\be
&&\mE\left(\sup\limits_{t\in[0,T]}|J_{11}(t)|\right)+\mE\left(\sup\limits_{t\in[0,T]}|J_{12}(t)|\right)+\mE\left(\sup\limits_{t\in[0,T]}|J_{13}(t)|\right)\no\\
&\leq& C\int_{0}^{T}\mE|Z^{\e,u_{\e}}(s)|^2\dif s+C\int_{0}^{T}\mE|X^{\e,\d,u_{\e}}_{s}-X^{\e,\d,u_{\e}}_{s(\triangle)}|^2\dif s\no\\
&&+C\int_{0}^{T}\mW_2^2(\sL_{X_s^{\e,\d}},\sL_{X_{s(\triangle)}^{\e,\d}})\dif s+C\int_{0}^{T}\mW_2^2(\sL_{X_s^{\e,\d}},D_{\bar{X}^{0}_{s}})\dif s\no\\
&&+C\int_{0}^{T}\mE|Y^{\e,\d,u_{\e}}_{s}-\hat{Y}^{\e,\d,u_{\e}}_{s}|^2\dif s\no\\
&\leq& C\int_{0}^{T}\mE|Z^{\e,u_{\e}}(s)|^2\dif s+C\left(\sup\limits_{s\in[0,T]}\mE\sup _{s \leqslant r \leqslant s+\triangle}|X_{r}^{\e,\d,u_\e}-X_{s}^{\e,\d,u_\e}|^{2}\right)\no\\
&&+C\left(\sup\limits_{s\in[0,T]}\mE\sup _{s \leqslant r \leqslant s+\triangle}|X_{r}^{\e,\d}-X_{s}^{\e,\d}|^{2}\right)+CT\mE\left(\sup\limits_{s\in[0,T]}|X_s^{\e,\d}-\bar{X}^{0}_{s}|^2\right)\no\\
&&+ CN\left(\frac{\d}{\e}\right).
\label{j111213}
\ee
And the H\"older inequality and the linear growth of $b_1, \bar{b}_1$ imply that
\be
\mE\left(\sup\limits_{t\in[0,T]}|J_{14}(t)|\right)
&\leq& C\left(\int_0^T(\mE|X^{\e,\d,u_{\e}}_{s}-X^{\e,\d,u_{\e}}_{s(\triangle)}|^2+\mE|\bar{X}^{u_\e}_{s}-\bar{X}^{u_\e}_{s(\triangle)}|^2)\dif s\right)^{1/2}\no\\
&&\times\left(\int_0^T\(1+\mE|X^{\e,\d,u_{\e}}_{s(\triangle)}|^2+\sL_{X_{s(\triangle)}^{\e,\d}}(|\cdot|^2)+\mE|\hat{Y}^{\e,\d,u_{\e}}_{s}|^2\)\dif s\right)^{1/2}\no\\
&\leq& C\Bigg(\left(\sup\limits_{s\in[0,T]}\mE\sup _{s \leqslant r \leqslant s+\triangle}|X_{r}^{\e,\d,u_\e}-X_{s}^{\e,\d,u_\e}|^{2}\right)\no\\
&&+\left(\sup\limits_{s\in[0,T]}\mE\sup _{s \leqslant r \leqslant s+\triangle}|\bar{X}^{u_\e}_r-\bar{X}^{u_\e}_s|^{2}\right)\Bigg)^{1/2}.
\label{j14}
\ee
For $J_{15}(t)$, we have that
\be
\mE\left(\sup\limits_{t\in[0,T]}|J_{15}(t)|\right)\leq C\((\frac{\d}{\triangle})^{1/2}+\triangle^{1/2}\).
\label{j15}
\ee
The proof of (\ref{j15}) is placed in the Appendix.

For $J_{2}$, by $(\mathbf{H}^{1'}_{b_{1}, \s_{1}})$ and the H\"older inequality, it holds that
\be
\mE\sup\limits_{t\in[0,T]}|J_{2}(t)|&\leq& 2\sqrt{L_{b_1,\s_1}}\mE\int_{0}^{T}|Z^{\e,u_{\e}}(s)|\(|Z^{\e,u_{\e}}(s)|+\mW_2(\sL_{X_s^{\e,\d}},D_{\bar{X}^{0}_{s}})\)|\dot{u}_{\e}(s)|\dif s\no\\
&\leq&C\mE\left(\int_{0}^{T}\(|Z^{\e,u_{\e}}(s)|+\mW_2(\sL_{X_s^{\e,\d}},D_{\bar{X}^{0}_{s}})\)|\dot{u}_{\e}(s)|\dif s\right)^2\no\\
&&+\frac{1}{4}\mE\left[\sup\limits_{t\in[0,T]}|Z^{\e,u_{\e}}(t)|^{2}\right]\no\\
&\leq&C\mE\left[\int_{0}^{T}|Z^{\e,u_{\e}}(s)|^{2}\dif s\right]+C\mE\left(\sup\limits_{t\in[0,T]}|X_t^{\e,\d}-\bar{X}^{0}_{t}|^2\right)\no\\
&&+\frac{1}{4}\mE\left[\sup\limits_{t\in[0,T]}|Z^{\e,u_{\e}}(t)|^{2}\right].
\label{j21}
\ee

For $J_{3}(t)$, from the Burkholder-Davis-Gundy inequality and the linear growth of $\s_1$, it follows that
\be
\mE\left(\sup\limits_{t\in[0,T]}|J_{3}(t)|\right)&\leq& 2\sqrt{\e}C\mE\left(\int_{0}^{T}|Z^{\e,u_{\e}}(s)|^2\|\s_1(X^{\e,\d,u_{\e}}_{s},\sL_{X_s^{\e,\d}}) \|^2 \dif s\right)^{1/2}\no\\
&\leq&\frac{1}{4}\mE\left[\sup\limits_{t\in[0,T]}|Z^{\e,u_{\e}}(t)|^{2}\right]+\sqrt{\e}C\mE\int_{0}^{T} \|\s_1(X^{\e,\d,u_{\e}}_{s},\sL_{X_s^{\e,\d}}) \|^2 \dif s\no\\
&\leq&\frac{1}{4}\mE\left[\sup\limits_{t\in[0,T]}|Z^{\e,u_{\e}}(t)|^{2}\right]+\sqrt{\e}C.
\label{j3}
\ee
For $J_{4}(t)$, by the linear growth of $\s_1$, we know
\be
\mE\left(\sup\limits_{t\in[0,T]}|J_{4}(t)|\right)\leq \e C\int_{0}^{T}(1+\mE|X^{\e,\d,u_{\e}}_{s}|^2+\mE|X^{\e,\d}_{s}|^2)\dif s \leq\e C.
\label{j4}
\ee

Combining (\ref{j111213})-(\ref{j4}) with (\ref{j1j2j3j4}), we can get
\ce
&&\mE\left(\sup\limits_{t\in[0,T]}|Z^{\e,u_{\e}}(t)|^{2}\right)\\
&\leq& C\int_{0}^{T}\mE\sup\limits_{r\in[0,s]}|Z^{\e,u_{\e}}(r)|^2\dif s+C\left(\sup\limits_{s\in[0,T]}\mE\sup _{s \leqslant r \leqslant s+\triangle}|X_{r}^{\e,\d,u_\e}-X_{s}^{\e,\d,u_\e}|^{2}\right)\no\\
&&+C\left(\sup\limits_{s\in[0,T]}\mE\sup _{s \leqslant r \leqslant s+\triangle}|X_{r}^{\e,\d}-X_{s}^{\e,\d}|^{2}\right)+CT\mE\left(\sup\limits_{s\in[0,T]}|X_s^{\e,\d}-\bar{X}^{0}_{s}|^2\right)\no\\
&&+ CN\left(\frac{\d}{\e}\right)+C\Bigg(\left(\sup\limits_{s\in[0,T]}\mE\sup _{s \leqslant r \leqslant s+\triangle}|X_{r}^{\e,\d,u_\e}-X_{s}^{\e,\d,u_\e}|^{2}\right)\no\\
&&+\left(\sup\limits_{s\in[0,T]}\mE\sup _{s \leqslant r \leqslant s+\triangle}|\bar{X}^{u_\e}_r-\bar{X}^{u_\e}_s|^{2}\right)\Bigg)^{1/2}+C\((\frac{\d}{\triangle})^{1/2}+\triangle^{1/2}\)+C(\sqrt{\e}+\e),
\de
which together with the Gronwall inequality implies that
\be
\mE\left(\sup\limits_{t\in[0,T]}|Z^{\e,u_{\e}}(t)|^{2}\right)\leq C\bigg[\Sigma(\e)+\frac{\d}{\e}+(\frac{\d}{\triangle})^{1/2}+\triangle^{1/2}+\sqrt{\e}+\e\bigg],
\label{zeue}
\ee
where
\ce
\Sigma(\e)&:=&\left(\sup\limits_{s\in[0,T]}\mE\sup _{s \leqslant r \leqslant s+\triangle}|X_{r}^{\e,\d,u_\e}-X_{s}^{\e,\d,u_\e}|^{2}\right)+\left(\sup\limits_{s\in[0,T]}\mE\sup _{s \leqslant r \leqslant s+\triangle}|X_{r}^{\e,\d}-X_{s}^{\e,\d}|^{2}\right)\\
&&+\Bigg(\left(\sup\limits_{s\in[0,T]}\mE\sup _{s \leqslant r \leqslant s+\triangle}|X_{r}^{\e,\d,u_\e}-X_{s}^{\e,\d,u_\e}|^{2}\right)+\left(\sup\limits_{s\in[0,T]}\mE\sup _{s \leqslant r \leqslant s+\triangle}|\bar{X}^{u_\e}_r-\bar{X}^{u_\e}_s|^{2}\right)\Bigg)^{1/2}\\
&&+\mE\left(\sup\limits_{s\in[0,T]}|X_s^{\e,\d}-\bar{X}^{0}_{s}|^2\right).
\de

{\bf Step 2.} We prove that for any $\eta>0$,
$$
\lim\limits_{\e\rightarrow0}\mP\left(\sup\limits_{t\in[0,T]}\left|\cG^{\e}\left(\sqrt{\e}W_\cdot+u_{\e}\right)(t)-\cG^{0}\left(u_\e\right)(t)\right|>\eta\right)=0.
$$

By the Chebyshev inequality, it holds that
\ce
&&\mP\left(\sup\limits_{t\in[0,T]}\left|\cG^{\e}\left(\sqrt{\e}W_\cdot+u_{\e}\right)(t)-\cG^{0}\left(u_\e\right)(t)\right|>\eta\right)\\
&=&\mP\left(\sup\limits_{t\in[0,T]}|Z^{\e,u_{\e}}(t)|>\eta\right)\leq\frac{1}{\eta^2}\mE\left(\sup\limits_{t\in[0,T]}|Z^{\e,u_{\e}}(t)|^{2}\right)\\
&\overset{(\ref{zeue})}{\leq}&C\frac{1}{\eta^2}\bigg[\Sigma(\e)+\frac{\d}{\e}+(\frac{\d}{\triangle})^{1/2}+\triangle^{1/2}+\sqrt{\e}+\e\bigg].
\de
Since $\lim\limits_{\e\rightarrow0}\d/\e=0$, $\d\rightarrow 0$ as $\e$ tends to $0$. Then we take $\triangle=\d^\g, 0<\g<1$ and have that $\triangle\rightarrow 0, \d/\triangle\rightarrow 0$, when $\e$ approximates to $0$. Hence, as $\e\rightarrow 0$, by (\ref{xegutse}) (\ref{barxuts}) and Theorem \ref{xbarxp}, it holds that 
$$
\lim\limits_{\e\rightarrow0}\mP\left(\sup\limits_{t\in[0,T]}\left|\cG^{\e}\left(\sqrt{\e}W_\cdot+u_{\e}\right)(t)-\cG^{0}\left(u_\e\right)(t)\right|>\eta\right)=0.
$$
The proof is complete.
\end{proof}

Now, it is the position to prove Theorem \ref{ldpmmsde}.

{\bf Proof of Theorem \ref{ldpmmsde}.}

By Theorem \ref{ldpbase}, to establish LDP, it is sufficient to verify that $\cG^0$ is measurable and Condition \ref{cond} $(i)$ and $(ii)$ hold. First of all, by Lemma \ref{auxilemm2}, we know that $\cG^0$ is measurable and Condition \ref{cond} $(i)$ is right. 
Then Lemma \ref{auxilemm3} implies that Condition \ref{cond} $(ii)$ is also right. The proof is complete.

\section{An example}\label{exam}

In this section, we explain our results by an example.

\bx
Consider the following slow-fast system of aggregation-diffusions equations on $\mR^n\times\mR^n$:
\be\left\{\begin{array}{l}
\dif X_{t}^{\e,\d}\in -\partial I_{\mathcal{O}}(X_{t}^{\e,\d})\dif t-\left[\nabla V_1(Y_{t}^{\e,\d})+\nabla V_2\ast\sL_{X_{t}^{\e,\d}}(X_{t}^{\e,\d})\right]\dif t+\sqrt{\e}\s_{1}\dif W^1_{t},\\
X_{0}^{\e,\d}=x_0\in\overline{\cD(\partial I_{\mathcal{O}})},\quad  0\leq t\leq T,\\
\dif Y_{t}^{\e,\d}\in -A_2(Y_{t}^{\e,\d})\dif t-\frac{1}{\d}\left[\nabla V_3(Y_{t}^{\e,\d})+\nabla V_4\ast\sL_{X_{t}^{\e,\d}}(X_{t}^{\e,\d})\right]\dif t+\frac{1}{\sqrt{\d}}\s_{2}\dif W^2_{t},\\
Y_{0}^{\e,\d}=y_0\in\overline{\cD(A_2)},\quad  0\leq t\leq T,
\end{array}
\right.
\label{Eqexldp}
\ee
where $\cO$ is a closed and convex domain in $\mR^n$ with ${\rm Int}(\cO)\neq\emptyset$, $V_i\in C^1(\mR^n)$ for $i=1,2,3,4$, $\ast$ denotes the convolution, $\s_1, \s_2$ are $n\times d_1, n\times d_2$ constant matrixes, respectively, and the rest of the setup is as in the system (\ref{Eq1}) with $n=m$.

Suppose that the derivative $\nabla V_1$ is Lipschitz continuous, the derivatives $\nabla V_2, \nabla V_4$ are bounded and Lipschitz continuous, and there exists a constant $\b>0$ such that for $y_1, y_2\in\mR^n$,
$$
\<y_1-y_2, \nabla V_3(y_1)-\nabla V_3(y_2)\>\leq -\b|y_1-y_2|^2.
$$
Then by Theorem \ref{xbarxpcr}, we know that for $0<\g<1$
\ce
\mE\(\sup_{0\leq t\leq T}|X_{t}^{\e,\d}-\bar{X}^0_{t}|^{2}\)\leq C(\d^{\g/2}+\d^\g+\d^{\frac{1}{2}(1-\g)}+\sqrt\e),
\de
where $(\bar{X}^{0},\bar{K}^{0})$ solves the following equation
\ce\left\{\begin{array}{l}
\dif \bar{X}^0_{t}\in -\partial I_{\mathcal{O}}(\bar{X}^0_{t})\dif t-\left[\int_{\cD(A_2)}\nabla V_1(y)\nu^{\bar{X}^0_{t},D_{\bar{X}^0_{t}}}(\dif y)+\nabla V_2(0)\right]\dif t,\\
\bar{X}^0_{0}=x_0\in\overline{\cD(\partial I_{\mathcal{O}})},
\end{array}
\right.
\de
and $\nu^{x,\mu}$ is the unique invariant probability measure of the following equation
\ce\left\{\begin{array}{l}
\dif Y_{t}^{x,\mu,y_0}\in -A_2(Y_{t}^{x,\mu,y_0})\dif t-\left[\nabla V_3(Y_{t}^{x,\mu,y_0})+\nabla V_4\ast\mu(x)\right]\dif t+\s_{2}\dif \hat{W}^2_{t},\\
Y_{0}^{x,\mu,y_0}=y_0\in\overline{\cD(A_2)}.
\end{array}
\right.
\de

Next, we consider the LDP for the system (\ref{Eqexldp}). Theorem \ref{ldpmmsde} implies that when $\lim\limits_{\e\rightarrow 0}\frac{\d}{\e}=0$, the family $\{X^{\e,\d},\epsilon\in(0,1)\}$ satisfies the LDP in $\mS:=C([0,T],\overline{\mathcal{D}(\partial I_{\mathcal{O}})})$ with the rate function given by
$$
I(\varsigma)=\frac{1}{2} \inf\limits_{h\in {\bf D}_{\varsigma}: \varsigma=\bar{X}^{h}}\|h\|_{\mH}^2,
$$
where $(\bar{X}^{h},\bar{K}^{h})$ solves the following equation
\ce\left\{\begin{array}{l}
\dif \bar{X}^h_{t}\in -A_1(\bar{X}^h_{t})\dif t-\left[\int_{\cD(A_2)}\nabla V_1(y)\nu^{\bar{X}^h_{t},D_{\bar{X}^0_{t}}}(\dif y)+\nabla V_2(\bar{X}^h_{t}-\bar{X}^0_{t})\right]\dif t+\s_{1}\pi_1\dot{h}(t)\dif t,\\
\bar{X}^h_{0}=x_0\in\overline{\cD(\partial I_{\mathcal{O}})}.
\end{array}
\right.
\de

Note that if the slow part of the system (\ref{Eqexldp}) doesn't depend on the fast part, i.e.
\ce
\left\{\begin{array}{l}
\dif X_{t}^{\e}\in -\partial I_{\mathcal{O}}(X_{t}^{\e})\dif t-\nabla V_2\ast\sL_{X_{t}^{\e}}(X_{t}^{\e})\dif t+\sqrt{\e}\s_{1}\dif W^1_{t},\\
X_{0}^{\e}=x_0\in\overline{\cD(\partial I_{\mathcal{O}})},\quad  0\leq t\leq T,
\end{array}
\right.
\de
the above equation falls into the class of equations in \cite{arrst}. There Adams et al. also studied the LDP under some similar assumptions.
\ex

\section{Appendix}\label{app}

In this section, we prove (\ref{gees}) and (\ref{j15}).

\subsection{Proof of (\ref{gees})} First of all, we justify that 

$(i)$ $\sup\limits_{\e\in (0,1)}\sup\limits_{t\in[0,T]}\left|g_\e(t)\right|<\infty$,

$(ii)$ $\{[0,T]\ni t\mapsto g_\e(t); \e\in (0,1)\}$ is equi-continuous.

For $0\leq s<t\leq T$, it holds that
\ce
|g_\e(s)-g_\e(t)|&\leq&\left|\int_{s}^{t}\sigma_1(\bar{X}^{h}_{r},D_{\bar{X}^0_{r}})(\pi_1h_{\e}(r)-\pi_1h(r))\dif r\right|\\
&\leq&\left(\int_{s}^{t}\|\sigma_1(\bar{X}^{h}_{r},D_{\bar{X}^0_{r}})\|^2\dif r\right)^{1/2}\left(\int_{s}^{t}|h_{\e}(r)-h(r)|^2\dif r\right)^{1/2}\\
&\leq&C\left(\int_{s}^{t}(1+|\bar{X}^{h}_{r}|^2+|\bar{X}^0_{r}|^2)\dif r\right)^{1/2}\times 2N^{1/2},
\de
where $({\bf H}_{b_1,\s_1}^{1'})$ is used. Letting $s=0$, by (\ref{barxub})  we have that
\ce
|g_\e(t)|\leq 2N^{1/2}C(1+|x_0|),
\de
where $C$ is independent of $\e$. So, $(i)$ holds.

For $(ii)$, noticing that 
\ce
|g_\e(s)-g_\e(t)|\leq 2N^{1/2}C(1+|x_0|)(t-s)^{1/2},
\de
we know that $(ii)$ holds.

Combining $(i)$ and $(ii)$, by the Ascoli-Arzel\'a lemma we obtain that $\{g_\e; \e\in (0,1)\}$ is relatively compact in $C([0,T],\mR^n)$.

Besides, note that 
\ce
\int_{0}^{t}\|\sigma_1(\bar{X}^{h}_{r},D_{\bar{X}^0_{r}})\|^2\dif r\leq C\int_{0}^{t}(1+|\bar{X}^{h}_{r}|^2+|\bar{X}^0_{r}|^2)\dif r<\infty. 
\de
Since $h_{\e}\rightarrow h$ in $\mathbf{D}_2^{N}$ as ${\e}\rightarrow0$, one get that for any $t\in[0,T]$
$$
\lim\limits_{{\e}\rightarrow0}g_\e(t)=0,
$$
which implies that
\ce
\lim\limits_{{\e}\rightarrow0}\sup\limits_{t\in[0,T]}\left|g_\e(t)\right|=0.
\de

\subsection{Proof of (\ref{j15})}

In this subsection, we prove (\ref{j15}). First of all, we collect some estimates for the frozen equation (\ref{Eq2}) which will be used in the sequel.

\bl\label{emb1}
 Suppose that $(\mathbf{H}^{1'}_{b_{1}, \s_{1}})$, $(\mathbf{H}_{A_{2}})$, $(\mathbf{H}^{2'}_{b_{2}, \s_{2}})$ and $(\mathbf{H}^{3}_{b_{2}, \s_{2}})$ hold. Then there exists a constant $C>0$ such that for any $t\in[0, T], x\in\overline{\cD(A_1)}, \mu\in\cP_2(\overline{\cD(A_1)}), y\in\overline{\cD(A_2)}$ 
\be
&&\hat{\mE}|Y_{t}^{x,\mu,y}|^{2}\leq|y|^{2}e^{-\frac{\a}{2} t}+C(1+|x|^{2}+\|\mu\|^2), \label{memu2}\\
&&\int_{\overline{\cD(A_2)}}|u|^2\nu^{x,\mu}(\dif u)\leq C(1+|x|^{2}+\|\mu\|^2),\\
&&|\hat{\mE} b_{1}(x,\mu,Y_{t}^{x,\mu,y})-\bar{b}_{1}(x,\mu)|^{2}\leq Ce^{-\a t}(1+|x|^{2}+\|\mu\|^2+|y|^{2}),
\label{meu2}
\ee
where $\hat{\mE}$ denotes the expectation under the probability measure $\hat{\mP}$.
\el

Since the proof of the above lemma is similar to that of \cite[Lemma 4.2 and 4.3]{q1}, we omit it.

{\bf Proof of (\ref{j15}).} For $J_{15}(t)$, it holds that
\ce
&&\mE\sup\limits_{t\in[0,T]}\left|J_{15}(t)\right|\\
&=&2\Bigg(\mE\sup_{0\leq t\leq T}
 \Big|\int_{0}^{[\frac{t}{\triangle}]\triangle}\<Z^{\e,u_{\e}}(s(\triangle)),b_{1}(X_{s(\triangle)}^{\e,\d,u_\e},\sL_{X_{s(\triangle)}^{\e,\d}},\hat{Y}_{s}^{\e,\d,u_\e})
 -\bar{b}_{1}(X_{s(\triangle)}^{\e,\d,u_\e},\sL_{X_{s(\triangle)}^{\e,\d}})\>\dif s\no\\
 &&\quad\quad\quad\quad+\int_{[\frac{t}{\triangle}]\triangle}^{t}\<Z^{\e,u_{\e}}(s(\triangle)),b_{1}(X_{s(\triangle)}^{\e,\d,u_\e},\sL_{X_{s(\triangle)}^{\e,\d}},\hat{Y}_{s}^{\e,\d,u_\e})
 -\bar{b}_{1}(X_{s(\triangle)}^{\e,\d,u_\e},\sL_{X_{s(\triangle)}^{\e,\d}})\>\dif s\Big|\Bigg)\no\\
&\leq&2\Bigg(\mE\sup_{0\leq t\leq T}
\Big|\int_{0}^{[\frac{t}{\triangle}]\triangle}\<Z^{\e,u_{\e}}(s(\triangle)),b_{1}(X_{s(\triangle)}^{\e,\d,u_\e},\sL_{X_{s(\triangle)}^{\e,\d}},\hat{Y}_{s}^{\e,\d,u_\e})
 -\bar{b}_{1}(X_{s(\triangle)}^{\e,\d,u_\e},\sL_{X_{s(\triangle)}^{\e,\d}})\>\dif s\Big|\Bigg)\no\\
&&+2\Bigg(\mE\sup_{0\leq t\leq T}
\Big|\int_{[\frac{t}{\triangle}]\triangle}^{t}\<Z^{\e,u_{\e}}(s(\triangle)),b_{1}(X_{s(\triangle)}^{\e,\d,u_\e},\sL_{X_{s(\triangle)}^{\e,\d}},\hat{Y}_{s}^{\e,\d,u_\e})
 -\bar{b}_{1}(X_{s(\triangle)}^{\e,\d,u_\e},\sL_{X_{s(\triangle)}^{\e,\d}})\>\dif s\Big|\Bigg)\no\\
 &=:&J_{151}+J_{152}.
\de

Next, we estimate $J_{151}$. Note that
\be
J_{151}
&=&2\Bigg(\mE\sup_{0\leq t\leq T}
\Big|\int_{0}^{[\frac{t}{\triangle}]\triangle}\<Z^{\e,u_{\e}}(s(\triangle)),\no\\
&&\qquad\qquad b_{1}(X_{s(\triangle)}^{\e,\d,u_\e},\sL_{X_{s(\triangle)}^{\e,\d}},\hat{Y}_{s}^{\e,\d,u_\e})
 -\bar{b}_{1}(X_{s(\triangle)}^{\e,\d,u_\e},\sL_{X_{s(\triangle)}^{\e,\d}})\>\dif s\Big|\Bigg)\no\\
&=&2\mE\Bigg(\sup_{0\leq t\leq T}
\Big|\sum\limits_{k=0}^{[\frac{t}{\triangle}]-1}\int_{k\triangle}^{(k+1)\triangle}\<Z^{\e,u_{\e}}(s(\triangle)),\no\\
&&\qquad\qquad b_{1}(X_{s(\triangle)}^{\e,\d,u_\e},\sL_{X_{s(\triangle)}^{\e,\d}},\hat{Y}_{s}^{\e,\d,u_\e})
 -\bar{b}_{1}(X_{s(\triangle)}^{\e,\d,u_\e},\sL_{X_{s(\triangle)}^{\e,\d}})\>\dif s\Big|\Bigg)\no\\
&\leq&2\mE\Bigg(\sup_{0\leq t\leq T}\sum\limits_{k=0}^{[\frac{t}{\triangle}]-1}\bigg{|}\int_{k\triangle}^{(k+1)\triangle}\<Z^{\e,u_{\e}}(k\triangle),\no\\
&&\qquad\qquad b_{1}(X_{k\triangle}^{\e,\d,u_\e},\sL_{X_{k\triangle}^{\e,\d}},\hat{Y}_{s}^{\e,\d,u_\e})-\bar{b}_{1}(X_{k\triangle}^{\e,\d,u_\e},\sL_{X_{k\triangle}^{\e,\d}})\>\dif s\bigg{|}\no\\
&\leq&2\sum\limits_{k=0}^{[\frac{T}{\triangle}]-1}\mE\Bigg(\Big|\int_{k\triangle}^{(k+1)\triangle}\<Z^{\e,u_{\e}}(k\triangle),\no\\
&&\qquad\qquad b_{1}(X_{k\triangle}^{\e,\d,u_\e},\sL_{X_{k\triangle}^{\e,\d}},\hat{Y}_{s}^{\e,\d,u_\e})
 -\bar{b}_{1}(X_{k\triangle}^{\e,\d,u_\e},\sL_{X_{k\triangle}^{\e,\d}})\>\dif s\Big|\Bigg)\no\\
&\leq&2[\frac{T}{\triangle}]\sup_{0\leq k\leq [\frac{T}{\triangle}]-1}
\mE\Bigg(\Big|\int_{k\triangle}^{(k+1)\triangle}\<Z^{\e,u_{\e}}(k\triangle),\no\\
&&\qquad\qquad b_{1}(X_{k\triangle}^{\e,\d,u_\e},\sL_{X_{k\triangle}^{\e,\d}},\hat{Y}_{s}^{\e,\d,u_\e})-\bar{b}_{1}(X_{k\triangle}^{\e,\d,u_\e},\sL_{X_{k\triangle}^{\e,\d}})\>\dif s\Big|\Bigg)\no\\
&\leq&2\d(\frac{T}{\triangle})\sup_{0\leq k\leq [\frac{T}{\triangle}]-1}\mE\Bigg(\Big|\<Z^{\e,u_{\e}}(k\triangle),\no\\
&&\qquad\qquad \int_{0}^{\triangle/\d}(b_{1}(X_{k\triangle}^{\e,\d,u_\e},\sL_{X_{k\triangle}^{\e,\d}},\hat{Y}_{\d s+k\triangle}^{\e,\d,u_\e})-\bar{b}_{1}(X_{k\triangle}^{\e,\d,u_\e},\sL_{X_{k\triangle}^{\e,\d}}))\dif s\>\Big|\Bigg)\no\\
&\leq&2\d(\frac{T}{\triangle})\sup_{0\leq k\leq [\frac{T}{\triangle}]-1}\mE|Z^{\e,u_{\e}}(k\triangle)|\no\\
&&\qquad\qquad \left|\int_{0}^{\triangle/\d}(b_{1}(X_{k\triangle}^{\e,\d,u_\e},\sL_{X_{k\triangle}^{\e,\d}},\hat{Y}_{\d s+k\triangle}^{\e,\d,u_\e})
 -\bar{b}_{1}(X_{k\triangle}^{\e,\d,u_\e},\sL_{X_{k\triangle}^{\e,\d}}))\dif s\right|\no\\
 &\leq&2\d(\frac{T}{\triangle})\sup_{0\leq k\leq [\frac{T}{\triangle}]-1}(\mE|Z^{\e,u_{\e}}(k\triangle)|^2)^{1/2}\no\\
 &&\times\left(\mE\left|\int_{0}^{\triangle/\d}(b_{1}(X_{k\triangle}^{\e,\d,u_\e},\sL_{X_{k\triangle}^{\e,\d}},\hat{Y}_{\d s+k\triangle}^{\e,\d,u_\e})
 -\bar{b}_{1}(X_{k\triangle}^{\e,\d,u_\e},\sL_{X_{k\triangle}^{\e,\d}}))\dif s\right|^2\right)^{1/2}.
\label{b41c}
\ee
Thus, we compute $\mE\left|\int_{0}^{\triangle/\d}(b_{1}(X_{k\triangle}^{\e,\d,u_\e},\sL_{X_{k\triangle}^{\e,\d}},\hat{Y}_{\d s+k\triangle}^{\e,\d,u_\e})-\bar{b}_{1}(X_{k\triangle}^{\e,\d,u_\e},\sL_{X_{k\triangle}^{\e,\d}}))\dif s\right|^2$. It is easy to see that
\ce
&&\mE\left|\int_{0}^{\triangle/\d}(b_{1}(X_{k\triangle}^{\e,\d,u_\e},\sL_{X_{k\triangle}^{\e,\d}},\hat{Y}_{\d s+k\triangle}^{\e,\d,u_\e})-\bar{b}_{1}(X_{k\triangle}^{\e,\d,u_\e},\sL_{X_{k\triangle}^{\e,\d}}))\dif s\right|^2\\
&=&2\mE\int_{0}^{\triangle/\d}\int_{r}^{\triangle/\d}\<b_{1}(X_{k\triangle}^{\e,\d,u_\e},\sL_{X_{k\triangle}^{\e,\d}},\hat{Y}_{\d s+k\triangle}^{\e,\d,u_\e})-\bar{b}_{1}(X_{k\triangle}^{\e,\d,u_\e},\sL_{X_{k\triangle}^{\e,\d}}),\\
&&\qquad b_{1}(X_{k\triangle}^{\e,\d,u_\e},\sL_{X_{k\triangle}^{\e,\d}},\hat{Y}_{\d r+k\triangle}^{\e,\d,u_\e})-\bar{b}_{1}(X_{k\triangle}^{\e,\d,u_\e},\sL_{X_{k\triangle}^{\e,\d}})\>\dif s\dif r.
\de
Therefore, set for $0<r<s\leq\triangle/\d$
\ce
\Psi(s,r)&:=&\mE\<b_{1}(X_{k\triangle}^{\e,\d,u_\e},\sL_{X_{k\triangle}^{\e,\d}},\hat{Y}_{\d s+k\triangle}^{\e,\d,u_\e})-\bar{b}_{1}(X_{k\triangle}^{\e,\d,u_\e},\sL_{X_{k\triangle}^{\e,\d}}),\\
&&\qquad b_{1}(X_{k\triangle}^{\e,\d,u_\e},\sL_{X_{k\triangle}^{\e,\d}},\hat{Y}_{\d r+k\triangle}^{\e,\d,u_\e})-\bar{b}_{1}(X_{k\triangle}^{\e,\d,u_\e},\sL_{X_{k\triangle}^{\e,\d}})\>,
\de
and we only investigate $\Psi(s,r)$.

For any $s>0$ and $\varpi\in L^2(\Omega,\mathscr{F}_s,\mP; \overline{\cD(A_1)}), \mu\in \cP_2(\overline{\cD(A_1)}), \varrho\in L^2(\Omega,\mathscr{F}_s,\mP; \overline{\cD(A_2)})$, we consider the following equation
\ce\left\{\begin{array}{l}
\dif \check{Y}_t^{\d, s, \varpi,\mu,\varrho}\in -A_2(\check{Y}_t^{\d, s, \varpi,\mu,\varrho})\dif t+\frac{1}{\d}b_2(\varpi,\mu,\check{Y}_t^{\d, s,\varpi,\mu,\varrho})\dif t+\frac{1}{\sqrt{\d}}\s_2(\varpi,\mu,\check{Y}_t^{\d, s, \varpi,\mu,\varrho})\dif W^2_t,\\
\check{Y}_s^{\d, s, \varpi,\mu,\varrho}=\varrho.
\end{array}
\right.
\de
Then it holds that
$$
\hat{Y}_t^{\e,\d,u_\e}=\check{Y}_t^{\d, k \triangle, X_{k \triangle}^{\e,\d,u_\e},\sL_{X_{k\triangle}^{\e,\d}}, \hat{Y}_{k \triangle}^{\e,\d,u_\e}}, \quad \hat{K}_t^{2,\e,\d}=\check{K}_t^{2,\d, k \triangle, X_{k \triangle}^{\e,\d,u_\e},\sL_{X_{k\triangle}^{\e,\d}}, \hat{Y}_{k \triangle}^{\e,\d,u_\e}}, \quad t \in[k \triangle,(k+1) \triangle],
$$
and furthermore
\ce
\Psi(s,r)&=&\mE\<b_{1}(X_{k\triangle}^{\e,\d,u_\e},\sL_{X_{k\triangle}^{\e,\d}},\check{Y}^{\d, k \triangle, X_{k \triangle}^{\e,\d,u_\e},\sL_{X_{k\triangle}^{\e,\d}}, \hat{Y}_{k \triangle}^{\e,\d,u_\e}}_{\d s+k\triangle})
 -\bar{b}_{1}(X_{k\triangle}^{\e,\d,u_\e},\sL_{X_{k\triangle}^{\e,\d}}),\\
 &&\qquad b_{1}(X_{k\triangle}^{\e,\d,u_\e},\sL_{X_{k\triangle}^{\e,\d}},\check{Y}^{\d, k \triangle, X_{k \triangle}^{\e,\d,u_\e},\sL_{X_{k\triangle}^{\e,\d}}, \hat{Y}_{k \triangle}^{\e,\d,u_\e}}_{\d r+k\triangle})
 -\bar{b}_{1}(X_{k\triangle}^{\e,\d,u_\e},\sL_{X_{k\triangle}^{\e,\d}})\>.
\de
Note that $X_{k \triangle}^{\e,\d,u_\e}, \hat{Y}_{k \triangle}^{\e,\d,u_\e}$ are $\sF_{k\triangle}$-measurable, and for any $x\in\overline{\cD(A_1)}, \mu\in \cP_2(\overline{\cD(A_1)}), y\in\overline{\cD(A_2)}$, $\check{Y}_t^{\d, k\triangle, x,\sL_{X_{k\triangle}^{\e,\d}},y}$ is independent of $\sF_{k\triangle}$. Thus, we have that
\ce
\Psi(s,r)
 &=&\mE\Bigg[\mE\Bigg[\<b_{1}(X_{k\triangle}^{\e,\d,u_\e},\sL_{X_{k\triangle}^{\e,\d}},\check{Y}^{\d, k \triangle, X_{k \triangle}^{\e,\d,u_\e},\sL_{X_{k\triangle}^{\e,\d}}, \hat{Y}_{k \triangle}^{\e,\d,u_\e}}_{\d s+k\triangle})
 -\bar{b}_{1}(X_{k\triangle}^{\e,\d,u_\e},\sL_{X_{k\triangle}^{\e,\d}}),\\
 &&b_{1}(X_{k\triangle}^{\e,\d,u_\e},\sL_{X_{k\triangle}^{\e,\d}},\check{Y}^{\d, k \triangle, X_{k \triangle}^{\e,\d,u_\e},\sL_{X_{k\triangle}^{\e,\d}}, \hat{Y}_{k \triangle}^{\e,\d,u_\e}}_{\d r+k\triangle})
 -\bar{b}_{1}(X_{k\triangle}^{\e,\d,u_\e},\sL_{X_{k\triangle}^{\e,\d}})\>\Bigg{|}\sF_{k\triangle}\Bigg]\Bigg]\\
&=&\mE\Bigg[\mE\Bigg[\<b_{1}(x,\sL_{X_{k\triangle}^{\e,\d}},\check{Y}^{\d, k \triangle,x,\sL_{X_{k\triangle}^{\e,\d}},y}_{\d s+k\triangle})
 -\bar{b}_{1}(x,\sL_{X_{k\triangle}^{\e,\d}}),\\
 &&\qquad\qquad b_{1}(x,\sL_{X_{k\triangle}^{\e,\d}},\check{Y}^{\d, k \triangle, x,\sL_{X_{k\triangle}^{\e,\d}},y}_{\d r+k\triangle}) -\bar{b}_{1}(x,\sL_{X_{k\triangle}^{\e,\d}})\>\Bigg]\Bigg{|}_{(x,y)=(X_{k\triangle}^{\e,\d,u_\e},\hat{Y}_{k \triangle}^{\e,\d,u_\e})}\Bigg].
\de

Here, we investigate $\check{Y}^{\d, k \triangle, x,\sL_{X_{k\triangle}^{\e,\d}},y}_{\d s+k\triangle}$. On one hand, it holds that
\ce
\check{Y}^{\d, k \triangle, x,\sL_{X_{k\triangle}^{\e,\d}},y}_{\d s+k\triangle}&=&y-\check{K}_{\d s+k\triangle}^{2,\d, k\triangle,x,\sL_{X_{k\triangle}^{\e,\d}},y}+\check{K}_{k\triangle}^{2,\d, k\triangle,x,\sL_{X_{k\triangle}^{\e,\d}},y}\\
&&+\frac{1}{\d} \int_{k\triangle}^{\d s+k\triangle}b_2(x,\sL_{X_{k\triangle}^{\e,\d}},\check{Y}^{\e, k \triangle, x,\sL_{X_{k\triangle}^{\e,\d}},y}_{r})\dif r\\
&&+\frac{1}{\sqrt{\d}} \int_{k\triangle}^{\d s+k\triangle} \s_2(x,\sL_{X_{k\triangle}^{\e,\d}},\check{Y}^{\d, k \triangle,x,\sL_{X_{k\triangle}^{\e,\d}},y}_{r})\dif W^2_r\\
&=&y-\check{K}_{\d s+k\triangle}^{2,\d, k\triangle, x,\sL_{X_{k\triangle}^{\e,\d}},y}+\check{K}_{k\triangle}^{2,\d, k\triangle, x,\sL_{X_{k\triangle}^{\e,\d}},y}\\
&&+\frac{1}{\d} \int_{0}^{\d s}b_2(x,\sL_{X_{k\triangle}^{\e,\d}},\check{Y}^{\d, k \triangle, x,\sL_{X_{k\triangle}^{\e,\d}},y}_{u+k\triangle})\dif u\\
&&+\frac{1}{\sqrt{\d}} \int_{0}^{\d s} \s_2(x,\sL_{X_{k\triangle}^{\e,\d}},\check{Y}^{\d, k \triangle,x,\sL_{X_{k\triangle}^{\e,\d}},y}_{u+k\triangle})\dif \tilde{W}^2_u\\
&=&y-\check{\tilde{\check{K}}}_{s}^{2,\d, k\triangle, x,\sL_{X_{k\triangle}^{\e,\d}},y}+\int_{0}^{s}b_2(x,\sL_{X_{k\triangle}^{\e,\d}},\check{Y}^{\d, k \triangle, x,\sL_{X_{k\triangle}^{\e,\d}},y}_{\d v+k\triangle})\dif v\\
&&+\int_{0}^{s} \s_2(x,\sL_{X_{k\triangle}^{\e,\d}},\check{Y}^{\d, k \triangle, x,\sL_{X_{k\triangle}^{\e,\d}},y}_{\d v+k\triangle})\dif \check{\tilde{W}}^2_v,
\de
where $\tilde{W}^2_u:=W^2_{u+k\triangle}-W^2_{k\triangle}$ and $\check{\tilde{W}}^2_v:=\frac{1}{\sqrt{\d}}\tilde{W}^2_{\d v}$ are two $m$-dimensional standard Brownian motions, and $\check{\tilde{\check{K}}}_{s}^{2,\d, k\triangle, x,\sL_{X_{k\triangle}^{\e,\d}},y}:=\check{K}_{\d s+k\triangle}^{2,\d, k\triangle, x,\sL_{X_{k\triangle}^{\e,\d}},y}-\check{K}_{k\triangle}^{2,\d, k\triangle, x,\sL_{X_{k\triangle}^{\e,\d}},y}$. On the other hand, note that the frozen equation (\ref{Eq2}) is written as
\ce
Y_{s}^{x,\sL_{X_{k\triangle}^{\e,\d}},y}=y-K^{2,x,\sL_{X_{k\triangle}^{\e,\d}},y}_s+\int_0^s b_{2}(x,\sL_{X_{k\triangle}^{\e,\d}},Y_{r}^{x,\sL_{X_{k\triangle}^{\e,\d}},y})\dif r+\int_0^s\s_{2}(x,\sL_{X_{k\triangle}^{\e,\d}},Y_{r}^{x,\sL_{X_{k\triangle}^{\e,\d}},y})\dif \hat W^2_{r}.
\de
Thus, for $s\in[0,\triangle/\d]$, $\check{Y}^{\d, k \triangle, x,\sL_{X_{k\triangle}^{\e,\d}},y}_{\d s+k\triangle}$ and $Y_{s}^{x,\sL_{X_{k\triangle}^{\e,\d}},y}$ have the same distribution, which implies that
\ce
&&\mE\Bigg[\<b_{1}(x,\sL_{X_{k\triangle}^{\e,\d}},\check{Y}^{\d, k \triangle, x,\sL_{X_{k\triangle}^{\e,\d}},y}_{\d s+k\triangle})-\bar{b}_{1}(x,\sL_{X_{k\triangle}^{\e,\d}}),b_{1}(x,\sL_{X_{k\triangle}^{\e,\d}},\check{Y}^{\d, k \triangle, x,\sL_{X_{k\triangle}^{\e,\d}},y}_{\d r+k\triangle})-\bar{b}_{1}(x,\sL_{X_{k\triangle}^{\e,\d}})\>\Bigg]\\
&=&\hat\mE\<b_{1}(x,\sL_{X_{k\triangle}^{\e,\d}},Y_{s}^{x,\sL_{X_{k\triangle}^{\e,\d}},y})-\bar{b}_{1}(x,\sL_{X_{k\triangle}^{\e,\d}}),b_{1}(x,\sL_{X_{k\triangle}^{\e,\d}},Y_{r}^{x,\sL_{X_{k\triangle}^{\e,\d}},y})-\bar{b}_{1}(x,\sL_{X_{k\triangle}^{\e,\d}})\>\\
&=&\hat\mE\Bigg[\hat\mE\Bigg[\Bigg<b_{1}(x,\sL_{X_{k\triangle}^{\e,\d}},Y_{s}^{x,\sL_{X_{k\triangle}^{\e,\d}},y})-\bar{b}_{1}(x,\sL_{X_{k\triangle}^{\e,\d}}), \\
&&\qquad\qquad b_{1}(x,\sL_{X_{k\triangle}^{\e,\d}},Y_{r}^{x,\sL_{X_{k\triangle}^{\e,\d}},y})-\bar{b}_{1}(x,\sL_{X_{k\triangle}^{\e,\d}})\Bigg>\Bigg{|}\sF_r^{\hat W^2}\Bigg]\Bigg]\\
&=&\hat\mE\Bigg[\Bigg<\hat\mE\left[b_{1}(x,\sL_{X_{k\triangle}^{\e,\d}},Y_{s}^{x,\sL_{X_{k\triangle}^{\e,\d}},y})\Bigg{|}\sF_r^{\hat W^2}\right]-\bar{b}_{1}(x,\sL_{X_{k\triangle}^{\e,\d}}),\\
&&\qquad\qquad b_{1}(x,\sL_{X_{k\triangle}^{\e,\d}},Y_{r}^{x,\sL_{X_{k\triangle}^{\e,\d}},y})-\bar{b}_{1}(x,\sL_{X_{k\triangle}^{\e,\d}})\Bigg>\Bigg]\\
&\leq&\left(\hat\mE\left|\hat\mE\left[b_{1}(x,\sL_{X_{k\triangle}^{\e,\d}},Y_{s}^{x,\sL_{X_{k\triangle}^{\e,\d}},y})\Bigg{|}\sF_r^{\hat W^2}\right]-\bar{b}_{1}(x,\sL_{X_{k\triangle}^{\e,\d}})\right|^2\right)^{1/2}\\
&&\qquad\qquad \left(\hat\mE|b_{1}(x,\sL_{X_{k\triangle}^{\e,\d}},Y_{r}^{x,\sL_{X_{k\triangle}^{\e,\d}},y})-\bar{b}_{1}(x,\sL_{X_{k\triangle}^{\e,\d}})|^2\right)^{1/2},
\de
where $\mathscr{F}_{r}^{\hat W^2}=\sigma\{\hat W^2_s,0\leq s \leq r\}\vee \cN$, and $\cN$ denotes  the collection of all $\mP$-zero sets. Moreover, based on (\ref{meu2}), we obtain that
\ce
&&\left(\hat\mE\left|\hat\mE\left[b_{1}(x,\sL_{X_{k\triangle}^{\e,\d}},Y_{s}^{x,\sL_{X_{k\triangle}^{\e,\d}},y})\Bigg{|}\sF_r^{\hat W^2}\right]-\bar{b}_{1}(x,\sL_{X_{k\triangle}^{\e,\d}})\right|^2\right)^{1/2}\\
&=&\left(\hat\mE\left|\hat\mE\left[b_{1}(x,\sL_{X_{k\triangle}^{\e,\d}},Y_{s-r}^{x,\sL_{X_{k\triangle}^{\e,\d}},\hat{y}})\right]\Bigg |_{\hat{y}=Y_{r}^{x,\sL_{X_{k\triangle}^{\e,\d}},y}}-\bar{b}_{1}(x,\sL_{X_{k\triangle}^{\e,\d}})\right|^2\right)^{1/2}\\
&\leq&\(Ce^{-\a(s-r)}(1+|x|^{2}+\sL_{X_{k\triangle}^{\e,\d}}(|\cdot|^2)+\mE|Y_{r}^{x,\sL_{X_{k\triangle}^{\e,\d}},y}|^{2})\)^{\frac{1}{2}}\\
&\leq&\(Ce^{-\a(s-r)}(1+|x|^{2}+\sL_{X_{k\triangle}^{\e,\d}}(|\cdot|^2)+|y|^{2})\)^{\frac{1}{2}},
\de
and
\ce
&&\left(\hat\mE\left|b_{1}(x,\sL_{X_{k\triangle}^{\e,\d}},Y_{r}^{x,\sL_{X_{k\triangle}^{\e,\d}},y})-\bar{b}_{1}(x,\sL_{X_{k\triangle}^{\e,\d}})\right|^2\right)^{1/2}\\
&=&\left(\hat\mE\left|b_{1}(x,\sL_{X_{k\triangle}^{\e,\d}},Y_{r}^{x,\sL_{X_{k\triangle}^{\e,\d}},y})-\int_{\overline{\cD(A_2)}}b_{1}(x,\sL_{X_{k\triangle}^{\e,\d}},u)\nu^{x,\sL_{X_{k\triangle}^{\e,\d}}}(\dif u)\right|^2\right)^{1/2}\\
&\leq&\left(\hat\mE\int_{\overline{\cD(A_2)}}\left|b_{1}(x,\sL_{X_{k\triangle}^{\e,\d}},Y_{r}^{x,\sL_{X_{k\triangle}^{\e,\d}},y})-b_{1}(x,\sL_{X_{k\triangle}^{\e,\d}},u)\right|^2\nu^{x,\sL_{X_{k\triangle}^{\e,\d}}}(\dif u)\right)^{1/2}\\
&\leq&C\left(\int_{\overline{\cD(A_2)}}\hat\mE|Y_{r}^{x,\sL_{X_{k\triangle}^{\e,\d}},y}-u|^2\nu^{x,\sL_{X_{k\triangle}^{\e,\d}}}(\dif u)\right)^{1/2}\\
&\leq&C\left(|y|^{2}e^{-\frac{\a}{2} t}+C(1+|x|^{2}+\sL_{X_{k\triangle}^{\e,\d}}(|\cdot|^2))\right)^{1/2}.
\de

Combining the above deduction, by (\ref{xeb}) (\ref{xeub}) (\ref{hatzub}) one can have that
\ce
\Psi(s,r)\leq Ce^{-\a(s-r)/2}.
\de
Inserting the above inequality in (\ref{b41c}), we get that
\be
J_{151}\leq C\d(\frac{T}{\triangle})\sup_{0\leq k\leq [\frac{T}{\triangle}]-1}\left(\int_{0}^{\frac{\triangle}{\d}}\int_{r}^{\frac{\triangle}{\d}}Ce^{-\a(s-r)/2}\dif s\dif r\right)^{1/2}\leq C(\frac{\d}{\triangle})^{1/2}.
\label{b4de}
\ee

Finally, we estimate $J_{152}$. By (\ref{b1line}) (\ref{xeub}) (\ref{hatzub}) and the H\"older inequality, one could get that
\ce
J_{152}
&\leq& 2\Bigg(\mE\sup_{0\leq t\leq T}
\int_{[\frac{t}{\triangle}]\triangle}^{t}|Z^{\e,u_{\e}}(s(\triangle))||b_{1}(X_{s(\triangle)}^{\e,\d,u_\e},\sL_{X_{s(\triangle)}^{\e,\d}},\hat{Y}_{s}^{\e,\d,u_\e})
 -\bar{b}_{1}(X_{s(\triangle)}^{\e,\d,u_\e},\sL_{X_{s(\triangle)}^{\e,\d}})|\dif s\Bigg)\\
 &\leq&2\triangle^{1/2}\left(\mE\sup_{0\leq t\leq T}\int_{[\frac{t}{\triangle}]\triangle}^{t}|b_{1}(X_{s(\triangle)}^{\e,\d,u_\e},\sL_{X_{s(\triangle)}^{\e,\d}},\hat{Y}_{s}^{\e,\d,u_\e})
 -\bar{b}_{1}(X_{s(\triangle)}^{\e,\d,u_\e},\sL_{X_{s(\triangle)}^{\e,\d}})|^2\dif s\right)^{1/2}\\ 
 &&\times\left(\mE\sup_{0\leq s\leq T}|Z^{\e,u_{\e}}(s)|^2\right)^{1/2}\\
 &\leq&2\triangle^{1/2}\left(\mE\int_{0}^{T}|b_{1}(X_{s(\triangle)}^{\e,\d,u_\e},\sL_{X_{s(\triangle)}^{\e,\d}},\hat{Y}_{s}^{\e,\d,u_\e})
 -\bar{b}_{1}(X_{s(\triangle)}^{\e,\d,u_\e},\sL_{X_{s(\triangle)}^{\e,\d}})|^2\dif s\right)^{1/2}\\ 
 &&\times\left(\mE\sup_{0\leq s\leq T}|Z^{\e,u_{\e}}(s)|^2\right)^{1/2}\\
&\leq& C\triangle^{1/2},
\de
which together with (\ref{b4de}) implies (\ref{j15}). The proof is complete.

\end{document}